\documentclass[11pt]{article}

\input{epsf.sty}

\usepackage{amsfonts,amssymb}

\newcommand{\qed}{\hbox{\rule{6pt}{6pt}}}

\newcommand{\Z}{\mathbb{Z}}

\newcommand{\Sym}{\mathbb{S}}

\newcommand{\Aut}{\mbox{Aut}}
\newcommand{\End}{\mbox{End}}
\newcommand{\tr}{\mbox{tr}}
\newcommand{\CC}{\mathbb{C}}
\newcommand{\BB}{\mathbb{B}}
\newcommand{\Sim}{\mathbb{S}}

\newtheorem{theorem}{Theorem}[section]

\newtheorem{lemma}[theorem]{Lemma}
\newtheorem{proposition}[theorem]{Proposition}

\newtheorem{definition}[theorem]{Definition}

\newtheorem{example}[theorem]{Example}

\newtheorem{remark}[theorem]{Remark}

\begin{document}

\title{
Cocycle Knot Invariants from Quandle Modules and 
Generalized Quandle Cohomology
}

\author{
J. Scott Carter\footnote{Supported in part by NSF Grant DMS \#9988107.}
\\University of South Alabama \\
Mobile, AL 36688 \\ carter@jaguar1.usouthal.edu
\and 
Mohamed Elhamdadi
\\ University of South Florida
\\ Tampa, FL 33620  \\ emohamed@math.usf.edu
\and 
Matias Gra\~{n}a\footnote{Supported in part by  CONICET and UBACyT X-066.} \\
Universidad de Buenos Aires \\
1428
Buenos Aires,  Argentina 
\\ matiasg@dm.uba.ar
\and 
Masahico Saito\footnote{Supported in part by NSF Grant DMS \#9988101.}
\\ University of South Florida
\\ Tampa, FL 33620  \\ saito@math.usf.edu
}

\maketitle

\begin{abstract} 
Three new knot invariants are defined using cocycles 
of the generalized quandle homology theory that was proposed by
Andruskiewitsch and Gra\~{n}a.
We specialize that theory to the case when there is a group action on
the coefficients.

First, quandle modules are used to generalize Burau representations
and Alexander modules for classical knots. Second, $2$-cocycles
valued in non-abelian groups are used in a way similar to Hopf algebra 
invariants of classical knots.
These invariants are shown to be of quantum type.
Third, cocycles with group actions on 
coefficient groups are used to define 
quandle cocycle 
invariants for both
classical knots and knotted surfaces. Concrete computational methods
are provided 
and used to prove non-invertibility for a large family of knotted surfaces.
In the classical case, the invariant 
can detect the chirality 
of 3-colorable knots in a number of cases.
\end{abstract}

\section{Introduction} \label{intro}

Quandle cohomology theory 
was developed \cite{CJKLS} to define invariants
of classical knots and knotted surfaces in state-sum form,
called  quandle cocycle (knot) invariants.
The quandle cohomology theory is a modification
of rack cohomology theory  which was defined in \cite{FRS},
and studied from different perspectives.
The cocycle knot invariants are analogous in their definitions
to the Dijkgraaf-Witten invariants
\cite{DW} of triangulated $3$-manifolds with finite gauge groups, 
but they 
use quandle knot colorings as spins and cocycles as 
Boltzmann weights.
 
Two types of topological applications of cocycle knot 
invariants have been established and are being investigated actively:
  non-invertibility  \cite{CJKLS,Satoh} and the
minimal triple point numbers  \cite{SatShi} of knotted surfaces.
A knot is {\it non-invertible} if it is not equivalent to itself with the 
orientation reversed 
while the 
the orientation of the space preserved.
In both applications, quandle cocycle invariants produced results 
that are not obtained by  other known methods; 
specifically, all known methods 
(e.g., \cite{Far75,Gor02*,Lev66,Rub83})
of proving non-invertibility of knotted surfaces
do not apply directly 
to non-spherical knots (except Kawauchi's \cite{Kawa86a,Kawa90a}
generalization of the Farber-Levine pairing, and it is not clear 
how this can be applied or 
how computations can be implemented for given examples), while 
quandle cocycle invariants can be applied  to orientable 
surfaces of any genus. Thus the cocycle invariants are the only 
method available in detecting non-invertibility 
that can be applied regardless of genus. 
Furthermore concrete 
methods 
to compute these have been implemented 
via computers.
The triple point numbers (the minimal number of triple points in 
projections)  have been  determined for some knotted surfaces
for the first time  
by using 
the  quandle cocycle invariants.

In this paper, we generalize the quandle cocycle invariants 
in three different directions, using generalizations of 
quandle homology theory provided by  Andruskiewitsch and Gra\~{n}a  
\cite{AG}. The original and the generalized quandle homology theories can be 
compared to group cohomology theories, with trivial and non-trivial 
group actions on coefficient groups, respectively. 
Thus the generalization of the homology theory and knot invariants are 
substantial and essential; the original case is only the very special case 
when the action is trivial. 
Examples in wreath product of groups are given in Section~\ref{extsec} 
after  preliminaries are provided in Section~\ref{prelimsec}.
Algebraic aspects of these examples are also studied.
The three directions of generalizations are as follows. 
First, the actions of quandle modules (defined below) are regarded as
generalizations of the Burau representation of braid groups. 
Thus we define invariants for classical knots, in Section~\ref{qmodinvsec},
 using such quandle modules, 
in a similar manner as Alexander modules are defined from 
Burau representations.   
Second,  quandle $2$-cocycles with non-abelian coefficients 
are used to define invariants for classical knots, 
defined in a similar manner as invariants defined from Hopf algebras
\cite{KauRad}, by sliding beads on knot diagrams, 
in Section~\ref{nonabsec}.
Generalizations of quandle cocycle invariants 
are given for classical knots in Section~\ref{cocyinvsec}
and for knotted surfaces in Section~\ref{knottedsfcesec},
respectively. Computational methods, examples, and applications are 
provided.
The invariant for the classical knots detect chirality 
of the $3$-colorable knots through $9$-crossings. 
As a main application of the invariant for knotted surfaces, 
we show that majority of $2k$-twist spun 
of $3$-colorable knots in the classical knot 
table up to $9$ crossings, as well as 
surfaces obtained from them by attaching trivial $1$-handles, 
are non-invertible.

\section{Preliminary} \label{prelimsec}

\subsection*{Quandles and knot colorings}

A {\it quandle}, $X$, is a set with a binary operation 
$(a, b) \mapsto a * b$
such that

(I) For any $a \in X$,
$a* a =a$.

(II) For any $a,b \in X$, there is a unique $c \in X$ such that 
$a= c*b$.

(III) 
For any $a,b,c \in X$, we have
$ (a*b)*c=(a*c)*(b*c). $

A {\it rack} is a set with a binary operation that satisfies 
(II) and (III).

Racks and quandles have been studied in, for example, 
\cite{Br88,FR,Joyce,K&P,Matveev}.
The axioms for a quandle correspond respectively to the 
Reidemeister moves of type I, II, and III 
(see
\cite{FR,K&P}, for example). 
Quandle structures
 have been found in areas other than knot theory, 
see \cite{AG} 
and \cite{Br88}.

In Axiom (II), the element  $c$
that is uniquely determined from   given $a, b \in X$ such that $a=c*b$,
is denoted by $c=a \bar{*} b$.
A function $f: X \rightarrow  Y$ between quandles
or racks  is a {\it homomorphism}
if $f(a \ast b) = f(a) * f(b)$ 
for any $a, b \in X$. 

The following are typical examples of quandles.  
A group $X=G$ with
$n$-fold conjugation
as the quandle operation: $a*b=b^{n} a b^{-n}$ or $a*b=b^{-n} a b^n$.
We denote by Conj$(G)$ the quandle defined for a group $G$ by 
$a*b=bab^{-1}$. 
Any subset of $G$ that is closed under such conjugation 
is also a quandle. 

Any $\Lambda (={\Z }[t, t^{-1}])$-module $M$  
is a quandle with
$a*b=ta+(1-t)b$, $a,b \in M$, 
that is 
called an {\it  Alexander  quandle}.
For a positive integer
$n$,
${\Z }_n[t, t^{-1}]/(h(t))$
is a quandle
for
a Laurent polynomial $h(t)$.
It 
is finite if the coefficients of the
highest and lowest degree terms
of $h$   are units in $\Z_n$.

Let $n$ be a positive integer, and 
for elements  $i, j \in \{ 0, 1, \ldots , n-1 \}$, define
$i\ast j \equiv 2j-i \pmod{n}$.
Then $\ast$ defines a quandle
structure  called the {\it dihedral quandle},
  $R_n$.
This set can be identified with  the
set of reflections of a regular $n$-gon
  with conjugation
as the quandle operation, but also is isomorphic to an Alexander quandle 
${\Z }_n[t, t^{-1}]/(t+1)$. 
As a set of reflections of the 
regular 
$n$-gon, 
 $R_n$ can be considered as a subquandle of ${\mbox{\rm Conj}}(\Sigma_n)$.

\begin{figure}
\begin{center}
\mbox{
\epsfxsize=2.5in
\epsfbox{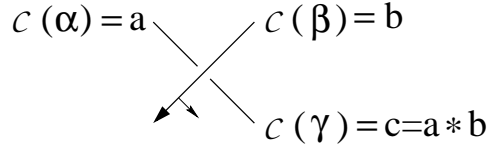} 
}
\end{center}
\caption{ Quandle relation at a crossing  }
\label{qcolor} 
\end{figure}

Let $X$ be a fixed quandle.
Let $K$ be a given oriented classical knot or link diagram,
and let ${\cal R}$ be the set of (over-)arcs. 
The normals are given in such a way that (tangent, 
normal) agrees with
the orientation of the plane, see Fig.~\ref{qcolor}. 
A (quandle) {\it coloring} ${\cal C}$ is a map 
${\cal C} : {\cal R} \rightarrow X$ such that at every crossing,
the relation depicted in Fig.~\ref{qcolor} holds. 
More specifically, let $\beta$ be the over-arc at a crossing,
and $\alpha$, $\gamma$ be under-arcs such that the normal of the over-arc
points from $\alpha$ to $\gamma$.
(In this case, $\alpha$ is called the {\it source arc} and  $\gamma$ 
is called the {\it target arc}.) 
Then it is required that ${\cal C}(\gamma)={\cal C}(\alpha)*{\cal C}(\beta)$.
The color ${\cal C}(\gamma)$ depends only on the choice 
of orientation of the over-arc; therefore this rule defines the coloring
 at both positive and negative crossings.
The  colors ${\cal C}(\alpha)$, ${\cal C}(\beta)$
are called {\it source} colors.

\subsection*{Quandle Modules}

We recall some information from \cite{AG}, 
but with notation changed to match our conventions. 

Let $X$ be a quandle. Let $\Omega (X)$ be the free ${\Z}$-algebra 
generated by 
$\eta_{x,y},$ $\tau_{x,y}$   for $x,y \in X$ such that 
$\eta_{x,y}$ is invertible for every $x,y \in X$.  
Define $\Z (X)$ to be the  quotient $\Z(X)= \Omega(X)/ R$ 
where $R$ is the ideal 
generated by 

\begin{enumerate}
\setlength{\itemsep}{-2pt}
\item \hfil $ \eta_{x*y,z}\eta_{x,y} -  \eta_{x*z,y*z}\eta_{x,z}$ \hfill \hfill
\item \hfil $ \eta_{x*y,z}\tau_{x,y} -  \tau_{x*z,y*z}\eta_{y,z}$ \hfill \hfill
\item \hfil $\tau_{x*y,z}- \eta_{x*z,y*z} \tau_{x,z}- 
 \tau_{x*z,y*z}\tau_{y,z}$ \hfill \hfill
\item \hfil $\tau_{x,x} + \eta_{x,x} -1$ \hfill \hfill
\end{enumerate}

The algebra  $\Z(X)$ thus defined is called the  {\it quandle algebra} 
over $X$. 
In $\Z(X)$, we define elements 
$\overline{\eta_{z,y}} =\eta^{-1}_{z\overline{*}y,y}$
and 
$\overline{\tau_{z,y}} = -\overline{\eta_{z,y}}\tau_{z\overline{*}y,y}$.
The convenience of such quantities will become apparent by examining type II 
moves.

A {\it representation} of $\Z(X)$ is an abelian group $G$ together with 
(1) a collection of automorphisms
$\eta_{x,y} \in {\mbox{\rm Aut}} (G)$, and (2) a collection of endomorphisms 
$\tau_{x,y}\in {\mbox{\rm End}} (G)$ such that the relations above hold. 
More precisely, there is an algebra homomorphism  
$\Z(X) \rightarrow {\mbox{\rm End}} (G)$, and we denote the 
image of the generators by the same symbols. Given a representation of 
$\Z(X)$ we say that $G$ is a $\Z(X)$-module, or {\it a quandle module}. 
The action of $\Z(X)$ on $G$ is written by the left action,  
and denoted by  $(\rho, g) \mapsto \rho g  (= \rho \cdot g = \rho (g) )$,
for  $g\in G$ and  $\rho \in End(G).$

\begin{example} {\bf \cite{AG}} \label{AGexample} {\rm 
Let $\Lambda = \Z[t,t^{-1}]$ denote the ring of Laurent polynomials. 
Then any  $\Lambda$-module $M$ is a $\Z(X)$-module for any quandle $X$,
by $\eta_{x,y} (a)=ta$ and $\tau_{x,y} (b) = (1-t) (b) $
for any $x,y \in X$.

The group  $G_X=\langle x \in X \ | \ x*y=yxy^{-1} \rangle$
is called the {\it enveloping group} \cite{AG}
(and the {\it associated group} in \cite{FR}). 
For any quandle $X$,
any $G_X$-module $M$  is a  $\Z(X)$-module by 
 $\eta_{x,y} (a)=y a$ and $\tau_{x,y} (b) = (1- x*y) (b) $,
where $x, y \in X$, $a, b \in M$.
} \end{example}

\begin{figure}[h]
\begin{center}
\mbox{
\epsfxsize=3in
\epsfbox{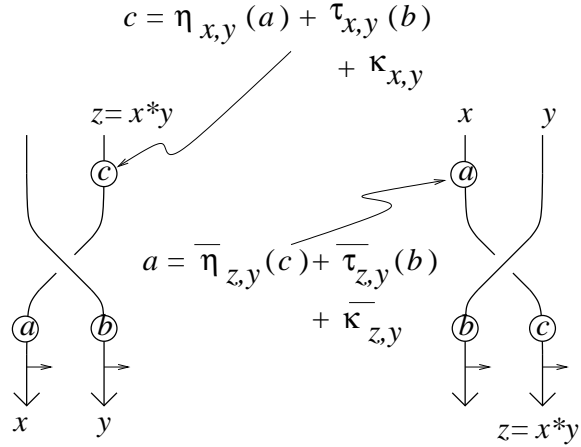} 
}
\end{center}
\caption{The geometric notation at a crossing}
\label{gencolor} 
\end{figure}

\begin{figure}[h]
\begin{center}
\mbox{
\epsfxsize=4in
\epsfbox{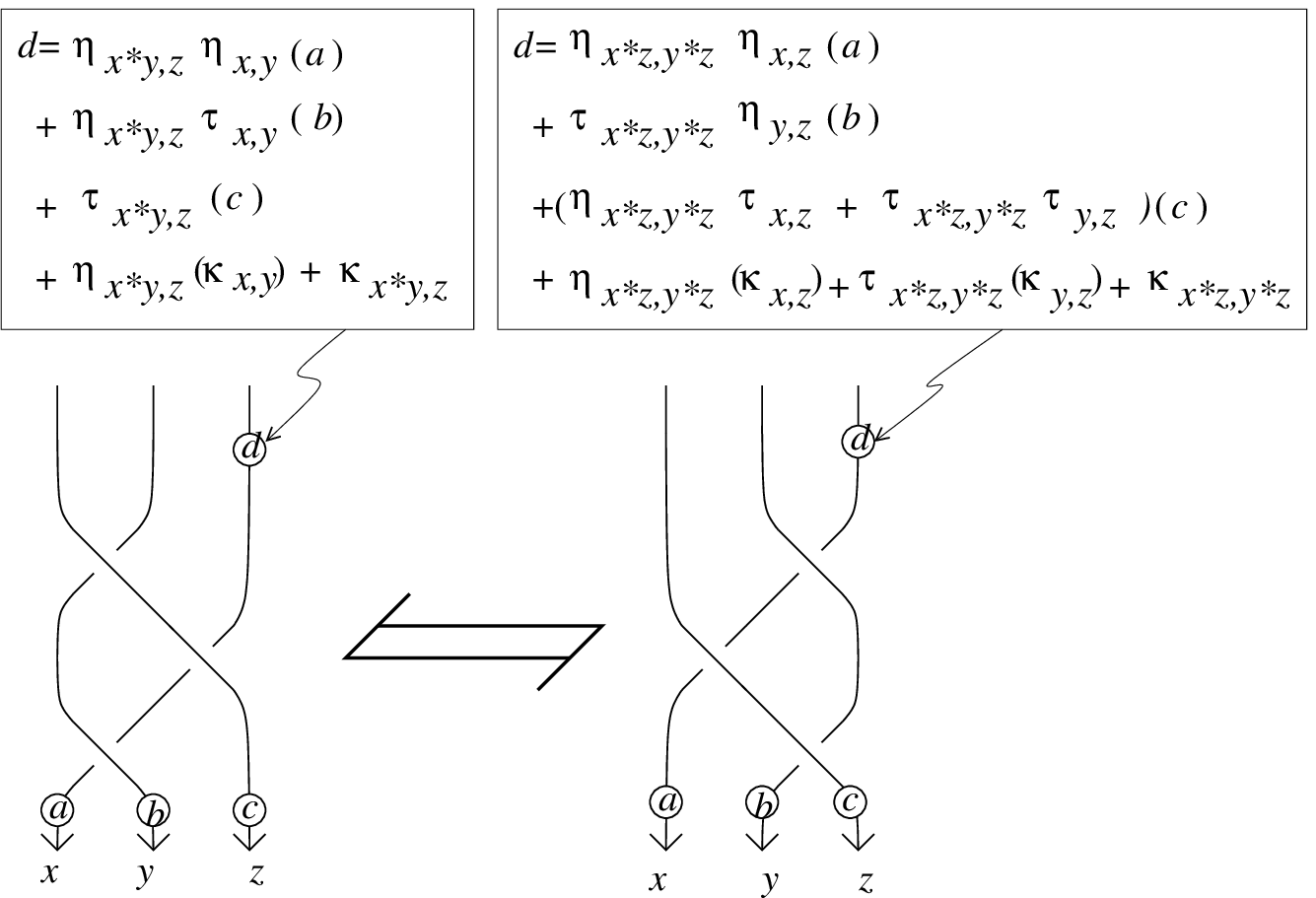} 
}
\end{center}
\caption{Reidemeister moves and the quandle algebra definition}
\label{move3} 
\end{figure}

We invite the reader to examine  
Figs.~\ref{gencolor} and \ref{move3} 
to see the geometric motivation for 
the quandle module axioms. 
For the time being, 
ignore the terms $\kappa_{x,y}$ in the figures. 
Detailed explanations of these figures will be given 
in  Section~\ref{qmodinvsec}. 

\subsection*{Generalized quandle homology theory}
Consider 
the free 
right 
$\Z ( X )$-module $C_n(X)=  \Z ( X ) X^n$ 
with basis $X^n$ (for $n=0$, $X^0$ is a singleton $\{ x_0 \}$, 
for a fixed element $x_0 \in X$). 
In \cite{AG},  boundary operators
 $\partial =\partial_n: C_{n+1}(X) \rightarrow C_n(X)$ 
are defined by 
\begin{eqnarray*} \partial (x_1, \ldots, x_{n+1} ) &=&
{\displaystyle (-1)^{n+1} \sum_{i=2}^{n+1} (-1)^i 
\eta_{ [x_1, \ldots, \widehat{x_i}, \ldots, x_{n+1} ], [x_i, \ldots, x_{n+1}] }
(x_1, \ldots, \widehat{x_i}, \ldots, x_{n+1} ) } \\
 & & \\
 & & - {\displaystyle (-1)^{n+1} \sum_{i=2}^{n+1} (-1)^i 
(x_1*x_i, \ldots, x_{i-1}*x_i, x_{i+1},  \ldots, x_{n+1} ) }  \\
 & & \\
 & & + (-1)^{n+1} 
 \tau_{[x_1, x_3, \ldots, x_{n+1}], [x_2, x_3, \ldots, x_{n+1}]}
(x_2, \ldots, x_{n+1} ) , \\[3mm]
\mbox{where} & & 
 [x_1, x_2, \ldots, x_n] = ( ( \cdots ( x_1 * x_2 ) * 
x_3  ) * \cdots ) * x_n
\end{eqnarray*}
for $n > 0$, and 
$\partial_1(x) = - \tau_{x \bar * x_0, x_0}$ for $n=0$.
The notational conventions  are slightly different from \cite{AG}.
In particular, the  $2$-cocycle 
condition for a $2$-cochain $\kappa_{x,y}$ 
in this homology theory 
is written as 
$$ \eta_{x*y, z} (\kappa_{x,y}) + \kappa_{x*y, z}
= \eta_{x*z, y*z} (\kappa_{x,z})  +  \tau_{x*z, y*z} (\kappa_{y,z})
 + \kappa_{x*z,y*z}, $$ 
for any $x,y,z \in X$.
We call this a {\it generalized $($rack$)$ $2$-cocycle condition}.
When $\kappa $ further satisfies $\kappa_{x,x}=0$ for any $x \in X$,
we call it a {\it generalized quandle $2$-cocycle.}

\subsection*{Dynamical cocycles}
Let $X$ be a quandle and $S$ be a non-empty set. 
Let $\alpha: X \times X \rightarrow 
\mbox{\rm Fun}(S \times S, S)=S^{S \times S}$ be a function,
so that for $x,y  \in X$ and $a, b \in S$ we have 
$\alpha_{x,y} (a,b) \in S$. 

Then it is checked by computations 
that $S \times X$ is a quandle by the operation
$(a, x)*(b, y)=(\alpha_{x, y}(a,b )  , x*y )$,
where $ x*y$ denotes the quandle 
product in $X$, 
if and only if $\alpha$ satisfies the following conditions:

\begin{enumerate}
\setlength{\itemsep}{-2pt}
\item $\quad\alpha_{x, x} (a,a)= a$ for all $x \in X$ and 
$a \in S$; 
\item $\quad\alpha_{x,y} (-,b): S \rightarrow S$ is a bijection for 
all $x, y \in X$ and for all $b \in S$;
\item $\quad 
\alpha_{x * y, z} ( \alpha_{x , y}(a, b), c)
=\alpha_{x * z, y*z}(\alpha_{x,z}(a,c), \alpha_{y,z}(b,c) )$
for all $x, y, z \in X$ and $a,b,c \in S$. 

\end{enumerate}

Such a function $\alpha$ is called a {\it dynamical quandle cocycle}
\cite{AG}.
The quandle constructed above is denoted by $S \times_{\alpha} X$, 
and is called the {\it extension} of $X$ by a dynamical cocycle $\alpha$.
The construction is general, as Andruskiewitsch and Gra\~{n}a
 show:

\begin{lemma}{\bf \cite{AG}} \label{AGlemma} 
Let $p: Y \rightarrow X$ be a surjective quandle homomorphism 
between finite quandles 
such that the cardinality of $p^{-1}(x)$ is a constant for all $x \in X$. 
Then $Y$ is isomorphic to an extension  $S \times_{\alpha} X$ of $X$ 
by some dynamical cocycle on 
the
set $S$ 
such that $|S|=|p^{-1}(x)|$.
\end{lemma}

\begin{example} {\bf \cite{AG} } 
{\rm
Let $G$ be  a $\Z(X)$-module for the quandle $X$,
and $\kappa$ be a generalized $2$-cocycle. 
For $a,b \in G$, let 
$$\alpha_{x,y}(a,b) =  \eta_{x,y}(a)  +  \tau_{x,y}(b)
+ \kappa_{x,y} .$$ 
Then it can be verified directly that $\alpha$ is a dynamical cocycle. 
In particular, even with $\kappa=0$, 
a $\Z(X)$-module structure on the abelian 
group $G$ defines a quandle structure
$G\times_\alpha X$. }\end{example}

\section{Group extensions and quandle modules} 
\label{extsec}
The purpose of this section is to provide examples of quandle modules
and cocycles from group extensions and group cocycles.  
Let  
$0\rightarrow N \stackrel{i}{\rightarrow} E 
\stackrel{\pi}{\rightarrow} H
\rightarrow 1$
be a 
short exact sequence of groups that expresses the 
group $E$ as a twisted
semi-direct product $E =N \rtimes_{\theta} H $ 
by a group $2$-cocycle $\theta$, 
where $N$ is an abelian group (see page $91$ of \cite{Brown}).
Thus
we have a set-theoretic section $s: H \rightarrow E $ 
that is normalized, in the sense that  $ s (1_H)= 1_{E} $, and  
the elements of $E$ 
can be written as 
pairs $(a, x)$  
 where $a\in N$ 
and $x \in H$,
by the bijection $(a, x) \mapsto i(a)s(x)$. 
We have 
$$s(x) s(y) = i (\theta (x, y) ) s(xy) , \quad x, y \in H, $$
and the multiplication 
rule in $E$  
is given by 
$(a, x) \cdot(b, y)=(a + x(b)+ \theta(x, y), xy)$, 
where $x(b)$ 
denotes the action of $H$ 
on $N$ that gives $E$  
the  structure of a twisted semi-direct product, $E=N \rtimes_{\theta} H$.
Recall here that the group $2$-cocycle condition is 
$$ \theta(x,y)+\theta(xy, z)=x \theta(y,z) + \theta(x, yz), \quad x,y,z \in H.$$

Let  $0\rightarrow N \stackrel{i}{\rightarrow} E 
\stackrel{\pi}{\rightarrow} H
\rightarrow 1$
be an exact sequence, and $E=N \rtimes_{\theta} H$ 
be the corresponding twisted semi-direct product.
Consider $E$ and $G$ as quandles, Conj$(E)$ and Conj$(H)$, 
respectively, where $a*b=bab^{-1}$. 
Lemma~\ref{AGlemma} implies that $E$ is  an extension of
$H$  by a dynamical cocycle $\alpha: H \times H \rightarrow N^{N \times N}$.

\begin{proposition}   \label{qmoduleprop}
The dynamical cocycle $\alpha$, 
in this case, 
is written by 
$\alpha_{x,y} (a,b ) =\eta_{x,y} (a) + \tau_{x,y}(b) + \kappa_{x,y}$
for any $x, y \in H$ and $a, b \in N$, where 
$\eta_{x,y} (a)=y  a $, 
$\tau_{x,y}(b)= (1-x*y) (b)$, and 
$$ \kappa_{x,y}=\theta(y,x) - yx \theta (y^{-1}, y) + \theta(yx, y^{-1}).$$
Similarly, if the quandle structure is defined by 
$r*s=s^{-1}rs$, then we obtain 
$\eta_{x,y} (a)=y^{-1}  a $, 
$\tau_{x,y}(b)= (y^{-1} x - y^{-1} ) (b)$ and
$$ \kappa_{x,y}= 
- \theta(y^{-1}, y) + y^{-1} \theta (x,y) + \theta(y^{-1}, xy) . $$ 
\end{proposition}
{\it Proof.\/} 
For $(b,y)\in E$, one has 
$$ (b, y)^{-1} = (-y^{-1} (b) - \theta(y^{-1}, y), y^{-1} ) 
=  (-y^{-1} (b) - y^{-1} \theta (y, y^{-1}), y^{-1} ) $$
(see page 92 of \cite{Brown}), and compute
\begin{eqnarray*} 
 (a, x)* (b, y) &=& 
(b, y) (a, x) (b, y)^{-1} \\
&= & ( b + y(a) - (yxy^{-1})(b) 
+ \theta(y,x)- yx \theta (y^{-1}, y) + \theta(yx, y^{-1}) , y x y^{-1}  )
\end{eqnarray*}
so that 
\begin{eqnarray*} 
\alpha_{x,y}(a,b)&=&\eta_{x,y}(a) + \tau_{x,y}(b) + \kappa_{x,y}\\
&=&
 y(a) +  (1-x*y) (b) + 
[\theta(y,x)- yx \theta (y^{-1}, y) + \theta(yx, y^{-1})] , \end{eqnarray*}
and we obtain the formulas.
Note that by expanding the terms $(a, x) (b, y)^{-1} $ first,
we obtain an equivalent formula 
$$\kappa_{x,y}=-yx\theta(y^{-1}, y)+ y \theta(x, y^{-1}) 
+ \theta(y, xy^{-1}) , $$
which follows from the group $2$-cocycle condition from the first 
formula, as well. 
The second case  is similar.
\qed

Thus the second item 
in Example~\ref{AGexample} 
occurs in semi-direct product of groups, when $\theta=0$ and 
hence $\kappa=0$.
Note also that the second case of Lemma~\ref{qmoduleprop}
agrees with the 
quandle action considered by Ohtsuki \cite{Ohtsuki}.

\begin{example} \label{wreathex} {\rm 
The wreath product of groups 
gives  
specific examples as follows.
Let 
$N=(\Z _q)^n$ for some $q\in \{ 0, 1, \ldots \}$.
(In case $q=0$, then $N$ is the direct product of the integers, and when 
$q=1$, then $N$ is trivial.)
The   symmetric group  $H=\Sigma_n$ acts 
on $N$  by  permutation of the factors
$\sigma(x_1,\ldots, x_n)= \sigma(\vec{x}) = 
(x_{\sigma^{-1} (1)}, \ldots, x_{\sigma^{-1} (n)})$, 
for $\sigma \in \Sigma_n$ and $\vec{x}=(x_j)_{j=1}^n \in N$. 
In this situation,
 $E$ is  called
a wreath product and denoted by 
$E= (\Z _q) \wr \Sigma_n$. 
In this case,
$\kappa=0$, and 
  the quandle module structure can be computed 
explicitly 
by matrices over $\Z _q$.
In \cite{CHNS},  such computations were used to obtain non-trivial 
colorings of some twist-spun knots by  dynamical extensions
of $R_n$. 
} \end{example}

Next we consider the $2$-cocycle $\kappa$ in terms of sections. 
In the group case, recall that the equality 
$s(x)s(y)=i(\theta (x,y)) s(xy)$ expresses
 the group $2$-cocycle $\theta$ as an obstruction to the section being 
a homomorphism. There is a similar  interpretation for quandle $2$-cocycles. 

\begin{lemma} \label{sectlemma}
The $2$-cocycle $\kappa$ 
in Proposition~\ref{qmoduleprop} 
satisfies
$s(x)*s(y)=i(\kappa_{x,y}) s(x*y)$. 
\end{lemma} 
{\it Proof.\/} 
One computes 
$$s(x)*s(y)=(0,x)*(0,y)=(\alpha_{x,y}(0,0), x*y)=(\kappa_{x,y}, x*y)
= i(\kappa_{x,y}) s(x*y)$$
as desired.
\qed

\begin{lemma} \label{kappathetalemma}
Let  $\kappa$ be as above. 
Then we have  $ \kappa_{x,y} = \theta(y, x) - \theta(yxy^{-1}, y ) .$
\end{lemma} 
{\it Proof.\/} 
{}From Lemma~\ref{sectlemma} we have 
\begin{eqnarray*} 
s(y) s(x) s(y)^{-1} &=& i\kappa_{x,y} s(yxy^{-1} ), \\
i\theta(y,x) s(yx) &=&i\kappa_{x,y} s(yxy^{-1} ) s(y) \\
 &=& i\kappa_{x,y}  i\theta(yxy^{-1}, y )s(yxy^{-1} y) ,
\end{eqnarray*}
and we obtain the formula, 
which is  simpler than that of Proposition~\ref{qmoduleprop}. \qed

\begin{lemma} \label{cobcoblemma}
Let $\kappa$ and $\theta$ be as above. 
If $\theta$ is a coboundary, then so is $\kappa$.
\end{lemma} 
{\it Proof.\/}
For a certain group $1$-cochain $\gamma$, 
$$ \theta(x,y)  
=\delta_G \gamma (x, y) 
=\gamma(xy)-\gamma(x)-x\gamma(y),$$
where $\delta_G$ (resp. $\delta_Q$) denotes the group
(resp. quandle) coboundary homomorphism. 
By Lemma~\ref{kappathetalemma},
$$\kappa_{x,y} 
= \gamma(yxy^{-1}) -y \gamma(x) - (1 - yxy^{-1})\gamma(y)
= \delta_Q \gamma(x,y)$$
as desired. \qed 

Thus  functors from group homology theories to quandle homology 
theories are expected.

Let   $E =N \rtimes_{\theta} H $ be as above, a twisted semi-direct product.
Let $X$ be a subquandle of Conj$(H)$,
and $\tilde{X}=\pi^{-1}(X)$, where $\pi: E \rightarrow H$ is the projection. 
Then $\tilde{X}$ is a subquandle of Conj$(E)$, 
and $\pi$ induces the quandle homomorphism 
$\pi:\tilde{X} \rightarrow X$, which is a dynamical extension.

\begin{example} \label{dihedexample} {\rm  
Let $X=R_n$ be the dihedral quandle of order $n$ (a positive integer),
which is a subquandle of Conj$(\Sigma_n)$.
Let $E=(\Z_q) \wr \Sigma_n$ be the wreath 
product as in Example~\ref{wreathex}, 
$$ 0 \rightarrow N=(Z_q)^n \rightarrow E=(\Z_q) \wr \Sigma_n
\stackrel{\pi}{\rightarrow} \Sigma_n \rightarrow 1. $$
Then $\tilde{X}=\pi^{-1}(X)$ is a subquandle of Conj$(E)$, and  
$\tilde{X}= (\Z_q)^n \times_{\alpha} X$ is a dynamical extension
of $X$. 
} \end{example}

As for the first cocycle 
groups, we have the following interpretation.
Let $X$ be a quandle, $A$ a $\Z(X)$-module. 
 Consider the dynamical extension $A \times_{\alpha} X$
of $X$ by $A$ with the dynamical cocycle 
$\alpha_{x,y}=\eta_{x,y}+\tau_{x,y}$ as before. 
For a given $1$-cochain $f \in C^1(X;A)=\mbox{\rm Hom}_{\Z(X)}(\Z(X)X,A)$, 
define a section
$\hat{f}: X \rightarrow A \times_{\alpha} X$ by 
$\hat{f}(x)=(f(x), x)$, which is 
indeed 
a section:
$\pi \circ \hat{f}={\rm id}_X$ for the projection 
$\pi: A \times_{\alpha} X \rightarrow X$.
\begin{lemma} \label{onecocylemma}
The section $\hat{f}$ is a quandle homomorphism if
and only if $f\in Z^1(X;A)$.
\end{lemma}
{\it Proof.\/}
The $1$-cocycle condition is written as
 $f(x*y)=\eta_{x,y}f(x)+\tau_{x,y}f(y)$
for any $x, y \in X$.
Thus we compute
\begin{eqnarray*}
\hat{f} (x*y) &=& (f(x*y), x*y) \\
\hat{f}(x) * \hat{f}(y) &=& 
(\alpha_{x,y}(f(x),f(y)), x*y) \ = \ 
(\eta_{x,y}f(x)+\tau_{x,y}f(y), x*y) \quad \mbox{\qed}.
\end{eqnarray*}

Alternatively, 
it can be stated that there is a one-to-one correspondence 
between  $Z^1(X;A)$ and the set of sections that are quandle homomorphisms.
 

\section{Knot invariants from quandle modules} \label{qmodinvsec}

\subsection*{Quandle modules and braids} \label{braidsec}

A braid word $w$ (of $k$-strings), or a $k$-braid word,  is 
a product of standard  generators $\sigma_1, \ldots, \sigma_{k-1}$ 
of the braid group ${\cal B}_k$ of $k$-strings and their inverses.
A braid word $w$ represents an element $[w]$ of the braid group 
${\cal B}_k$.
Geometrically, $w$ is represented by a diagram in a rectangular box
with  $k$ end points  at the top, and $k$ end points at the bottom,
where the strings   go down monotonically. 
Each generator or its inverse is represented by a crossing in a diagram.
We use the same letter $w$ for a choice of such a diagram. 
Let $\hat{w}$  denote the closure of 
the 
diagram 
$w$. 
Quandle colorings of 
$w$ 
are  defined in 
exactly the same manner as in the case of knots.
However, the quandle elements at the top and the bottom of a diagram of 
$w$ do not necessarily coincide. 
For the closure $\hat{w}$, the 
 quandle elements at the top and the bottom of a diagram of 
$w$ coincide, when we consider a coloring of a link  $\hat{w}$.

\begin{figure}[h]
\begin{center}
\mbox{
\epsfxsize=1in
\epsfbox{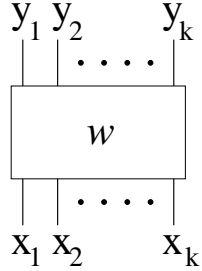} 
}
\end{center}
\caption{A quandle coloring of a braid word $w$}
\label{bcolors} 
\end{figure}

Let $X$ be a quandle.
Let $\gamma_1, \ldots, \gamma_k$ be the bottom arcs of $w$.
For a given vector $\vec{x}=(x_1, \ldots, x_k) \in X^k$,
assign these elements $ x_1, \ldots, x_k$ on $\gamma_1, \ldots, \gamma_k$ 
 as their colors, respectively. Then from the definition, 
a coloring ${\cal C}$ of $w$ by $X$ is uniquely determined
such that  
${\cal C}(\gamma_i)=x_i$, $i=1, \ldots, k$.
We call such a coloring ${\cal C} $ the {\it  coloring induced from} $\vec{x}$.
Let $\delta_1, \ldots, \delta_k$ be the arcs at the top.
Let 
$\vec{y}=(y_1, \ldots, y_k)=({\cal C}(\delta_1), \ldots, {\cal C}(\delta_k) )\in X^k$ 
be the colors assigned to the top arcs, that are uniquely determined from 
$\vec{x}$.
Denote
 this situation by a left action, 
$\vec{y} =  w \cdot \vec{x}$. The colors $\vec{x}$ and $\vec{y}$ are 
called {\it bottom} and {\it top colors}
or {\it color vectors}, 
respectively.
See Fig.~\ref{bcolors}. 

Let $X$ be a quandle and $G$ be a quandle module. 
For a  dynamical cocycle 
$\alpha= \eta + \tau + \kappa$ --- which acts on $(a,b) \in G^2$ by
$ \alpha_{x,y}(a,b) =\eta_{x,y} (a )  +  \tau_{x,y}(b) + \kappa_{x,y}$
for any $(x,y) \in X^2$ --- 
let 
$\tilde{X}=G \times_{\alpha} X$ be the dynamical extension. 
If $\vec{r} =( (a_1, x_1), \ldots, (a_k, x_k)) $
and $\vec{s} =( (b_1, y_1), \ldots, (b_k, y_k)) \in \tilde{X}^k$
are bottom and top colors of $w \in {\cal B}_k$ by $\tilde{X}$,
respectively, then we write this situation by 
$\vec{b} =  M(w, \vec{x}) \cdot \vec{a}$, 
where $\vec{a}=(a_1, \ldots, a_k)$, 
$\vec{b}=(b_1, \ldots, b_k) \in G^k$. 
Thus $M(w, \vec{x})$ represents a map 
$M(w, \vec{x}): G^k \rightarrow G^k$.

\begin{lemma} \label{indeplemma}
If $[w]=[w'] \in {\cal B}_k$, then 
$M(w, \vec{x})=M(w', \vec{x}) : G^k \rightarrow G^k$. 
\end{lemma}
{\it Proof.\/}
The invariance under the braid relations are checked from the definitions.
In particular,  
\begin{eqnarray*} 
\lefteqn{  M( \sigma_1 \sigma_2 \sigma_1 , (x,y,z) )(a,b,c)  } \\
 &= & (c,\  \eta_{y,z}( b) +   \tau_{y,z}(c) + \kappa_{y,z} , \
   \eta_{x*y, z} \eta_{x,y}(a) +  \eta_{x*y, z} \tau_{x,y}(b)
+  \tau_{x*y, z}(c) +  \eta_{x*y, z}( \kappa_{x,y}) + \kappa_{x*y, z} ), \\
\lefteqn{ M( \sigma_2 \sigma_1 \sigma_2 , (x,y,z) ) (a,b,c)  } \\
 &= & (c, \   \eta_{y,z}(b) +   \tau_{y,z}(c) + \kappa_{y,z} , \\
 & &  \eta_{x*z, y*z} \eta_{x,z}(a) +  \tau_{x*z, y*z} \eta_{y,z}(b)
+  ( \eta_{x*z, y*z} \tau_{x,z} +  \tau_{x*z, y*z}\tau_{y,z} )(c) \\ & &  
+ \eta_{x*z, y*z}( \kappa_{x,z}) + \tau_{x*z, y*z} ( \kappa_{y,z}) +
\kappa_{x*z, y*z} )
\end{eqnarray*}
and the equality follows from the quandle module conditions and 
the generalized $2$-cocycle condition.
\qed

This is not a braid group representation on $G^k$, 
as it depends on the color of $w$ by $X$.
However, in the case in which the coloring by $X$ is trivial
 so 
$x_1=x_2=\cdots = x_k$, and $\kappa=0$, then
it is a braid group representation.  We call the map 
$M(-, \vec{x}): {\cal B}_k \rightarrow \mbox{Map}(G^k, G^k)$ 
a {\it colored representation}.

For a standard braid generator, this situation is diagrammatically 
represented as depicted in Fig.~\ref{gencolor},
the left figure for a negative crossing, and the right one for positive. 
In the  calculations given in this section,
 the left figure represents the braid generator $\sigma_j$, and the 
right represents the inverse. 
In the figure, the colors by $X$ are assigned to arcs. 
Elements of $G$ are put in small circles on arcs. 
We 
imagine these circles sliding 
up through a crossing, at which 
the dynamical cocycle $\alpha$ acts and changes the elements
when a circled elements goes under a crossing. 
Going over a crossing does not change the element in a circle.
{}From type II Reidemeister moves, the definition
of $\overline{\eta}$ and $\overline{\tau}$ is recovered. 
Figure~\ref{move3} shows that the quandle module 
conditions correspond to the type III  move with this diagrammatic convention.

\subsection*{ Module invariants} 
\label{qmodulesec} 

Let $w$  be a $k$-braid word, and denote by $\hat{w}$ the 
closure of $w$.
Let $X$ be a quandle and $G$ be a quandle module.
Let $\alpha= \eta + \tau $ 
be a dynamical cocycle, which acts on $(a,b) \in G^2$ by
$ \alpha_{x,y}(a,b) =\eta_{x,y} (a )  +  \tau_{x,y}(b)$ 
for any $(x,y) \in X^2$. 

\begin{theorem}
Let $L$ be a link represented as a closed braid $\hat{w}$,
where $w$ is a $k$-braid word,   and 
${ \mbox{Col}_X(L) }$ be the set of colorings of $L$ by a quandle $X$.
For ${\cal C} \in  \mbox{Col}_X(L)$, 
let $\vec{x}$ be the color vector of bottom strings of 
$w$ that is the restriction of ${\cal C}$. 
Then the family  
$$\tilde{\Phi}(X, \alpha\ ; L)=\{ 
 G^k 
/\mbox{\rm Im} (M(w,\vec{x})-I)  \}_{{\cal C} \in  {\rm Col}_X(L)}$$
of isomorphism classes of modules presented by 
the maps $(M(w,\vec{x})-I)$, 
where $I$ denotes the identity,  is independent of choice of $w$
that represents $L$ as its closed braid,
and thus defines a link invariant.
\end{theorem}
{\it Proof.\/}  
 By the Lemma~\ref{indeplemma}, 
$M(w)=M(w, \vec{x})$ does not depend on the choice of 
a braid word. We use Markov's theorem to prove the statement.
First note that the set of colorings remains unchanged by a 
stabilization, in the sense that  there is a natural bijection
(the colorings on bottom  strings $(x_1, \ldots, x_k)$ of
a braid word $w$ extend uniquely
to the coloring $(x_1, \ldots, x_k, x_{k})$
of $w\sigma_{k}^{\pm 1}$, a stabilization of $w$).
There is a bijection of colorings between conjugates, as well.
Hence it is sufficient to prove that,  for a given coloring,
the isomorphism class of the module defined in the statement
remains unchanged by 
conjugation and stabilization for a natural induced coloring.
The  invariance under conjugation is seen from  the fact that 
conjugation by a braid word 
induces a 
conjugation 
by a matrix, and the module 
is isomorphic under conjugate presentation matrices.
So we investigate the stabilization.

We represent maps of $G^k$ by $k$ by $k$ matrices whose 
entries represent maps  of $G$.  
The braid generator $\sigma_k$
in the stabilization  $w\sigma_{k}^{\pm 1}$ 
 is represented by   the 
matrix 
$M(\sigma_k)=I_{k-1} 
\oplus 
\left[ \begin{array}{cc} O & I\\ W & I - W \end{array} \right]$,
 where $I$ denotes the identity map of $G$, 
and $I_k$ denotes the identity map on $G^k$.  
This is because the $k\/$th and $(k+1)\/$st strings receive 
the same color. 
The block matrix 
$\left[ \begin{array}{cc} O & I \\ W &  I -W \end{array} \right]$,
in particular, represents the map 
$$(a,b) \mapsto (b, \alpha_{x,x}(a,b))
=(b, \eta_{x,x}(a) + \tau_{x,x}(b) )$$
where $x=x_k$, so that 
 $W$ corresponds to the action by $\eta$,
and hence $W$ is an isomorphism of $G$. 
The fact that $\tau_{x,x}$
corresponds to the matrix $I - W$  
follows from the condition 
$\eta_{x,x}+\tau_{x,x}=1$ in a quandle algebra.

Express the matrix $M(w)$, where $w$ is regarded as 
a $(k+1)$-braid word after stabilization, 
 though originally a $k$-braid word,
 by a matrix 
$\left[ \begin{array}{cccc} M_{11} & \cdots &  M_{1k} &O\\ 
\vdots & & \vdots &  \vdots\\
M_{k1} &  \cdots & M_{kk} & O
\\ O &  \cdots & O & I \end{array} \right]$
where $M_{ij}$, $1 \leq i, j \leq k$, are maps of $G$, 
and $M(w)$ was written as  a $(k+1) \times (k+1)$ matrix of
these maps, 
$I$ denotes the 
identity map on $G$. 
Then $M(w \sigma_k)$ is represented by 
the matrix $$ M(w) M(\sigma_k)=
\left[ \begin{array}{ccccc} 
M_{11} &  \cdots  & M_{1\ (k-1)} &O & M_{1 k} \\ 
\vdots &  & \vdots & \vdots & \vdots \\
M_{k1} &  \cdots & M_{k\ (k-1)} & O & M_{k k} \\
O& \cdots  & O& W & I - W 
 \end{array} \right] .  $$
Hence 
$$ M(w) M(\sigma_k) - I_{k} = 
\left[ \begin{array}{ccccc} 
M_{11} - I  &  \cdots  & M_{1\ (k-1)} &O & M_{1 k} \\ 
\vdots &  & \vdots & \vdots & \vdots \\
M_{(k-1)\ 1} & \cdots &  M_{(k-1)\ (k-1)} - I & O &  M_{(k-1)\ k} \\
M_{k1} &  \cdots & M_{k\ (k-1)} & -I  & M_{k k} \\
O& \cdots  & O& W &  - W 
 \end{array} \right] ,
  $$
which is column reduced to 
$
\left[ \begin{array}{ccccc} 
M_{11} - I  &  \cdots  & M_{1\ (k-1)} &O & M_{1 k} \\ 
\vdots &  & \vdots & \vdots & \vdots \\
M_{(k-1)\ 1} & \cdots &  M_{(k-1)\ (k-1)} - I & O &  M_{(k-1)\ k} \\
M_{k1} &  \cdots & M_{k\ (k-1)} & -I  & M_{k k} -I \\
O& \cdots  & O& W &  O 
 \end{array} \right] 
  $
and is further row reduced to 
$
\left[ \begin{array}{ccccc} 
M_{11} - I  &  \cdots  & M_{1\ (k-1)} &O & M_{1 k} \\ 
\vdots &  & \vdots & \vdots & \vdots \\
M_{(k-1)\ 1} & \cdots &  M_{(k-1)\ (k-1)} - I & O &  M_{(k-1)\ k} \\
M_{k1} &  \cdots & M_{k\ (k-1)} & O  & M_{k k} -I \\
O& \cdots  & O& W &  O 
 \end{array} \right] 
  $
that represents the isomorphic module as $M(w)$ does, 
since $W$ is an isomorphism. 
The invariance under stabilization by $\sigma_k^{-1}$ follows similarly.
\qed

This theorem implies that the following is well-defined.

\begin{definition} {\rm 
The family of modules
$\tilde{\Phi}(X, \alpha \ ; L)
=\{ G^k 
/\mbox{\rm Im}(M(w,\vec{x})-I)\}_{{\cal C} \in {\rm Col}_X(L) } $
is called the {\em quandle module invariant}.
} \end{definition}

Specific examples can be constructed from wreath products.
Let $X$ be a subquandle of Conj$(\Sigma_n)$, 
and $E$ be a wreath product extension by $N=(\Z_q)^n$
where $q, n$ are positive integers, so that  
$$ 0 \rightarrow N=(\Z_q)^n \stackrel{i}{\rightarrow} E 
\stackrel{\pi}\rightarrow \Sigma_n \rightarrow 1, $$
and $E=(\Z_q )^n \rtimes \Sigma_n$.
Let  
$\tilde{X}=\pi^{-1}(X)$, as before, 
and let $\alpha$ be the corresponding dynamical cocycle,
so that $\tilde{X}=(\Z_q)^n \times_{\alpha} X  $.

\begin{example}{\rm 
Let $X=R_n$ denote the $n$-element dihedral quandle, regarded as a subset
of $\Sigma_n$. The action on $\Z^n$ is by permutations. 
Then the quandle module invariant is defined as above, with $q=0.$
We ran programs in {\it Maple}, 
 {\it Mathematica}, and/or {\it C}, 
 independently, 
to compute the quandle module invariants, for $n=3,5,7,11,$ and $13$,
 through the nine-crossing knots.
There are $84$ such prime knots in the table. 
We used Jones's~\cite{Jones} 
table in 
 which knots are 
given in braid form.  
For $n=3$, $5$, $7$, $11$, and $13$, 
there are $33$, $17$, $10$, and $7$ knots that are non-trivially colored 
by $R_n$, respectively.  
We summarize our data in Table~\ref{qmoduletable3}
and  ~\ref{qmoduletable5}
for $R_3$ and $R_5$, respectively, 
below.
For each knot we list the torsion subgroups, 
the rank of the  free part, and whether the coloring is trivial or not.
 Thus the entry in Table A-1 
for knot $8_{15}$ indicates that there are 
three trivial colorings of $8_{15}$ which give the invariant 
$\Z_{33}\oplus \Z^3$, and 
six colorings give the invariant $\Z_2 \oplus \Z^4$. 
For colorings by 
$R_7$, $R_{11}$ and $R_{13}$,  
the results are summarized as follows, where $D$ denotes the determinant
of a given knot.
\begin{itemize}
\setlength{\itemsep}{-2pt}
\item
All $7$-colorable knots up to $9$-crossings, except $9_{41}$,
have the module invariant $(\Z_D)^3 \oplus \Z^7$, for $7$ trivial 
colorings, and $\Z^{10}$ for $42$ non-trivial colorings.
For $9_{41}$, it is $(\Z_D)^6 \oplus \Z^7$  for $7$ trivial 
colorings, and $\Z^{10}$ for $336$ non-trivial colorings.

\item
All $11$-colorable knots up to $9$-crossings
have the module invariant $(\Z_D)^5 \oplus \Z^{11}$, for $11$ trivial 
colorings, and $\Z^{16}$ for $110$ non-trivial colorings.

\item
All $13$-colorable knots up to $9$-crossings
have the module invariant $(\Z_D)^6 \oplus \Z^{13}$, for $13$ trivial 
colorings, and $\Z^{19}$ for $156$ non-trivial colorings. 
\end{itemize}

We expect that the monotony of values for $7$, $11$, 
and $13$ is due to the fact that the knots considered 
are all of relatively small crossing numbers and/or bridge numbers.

} \end{example}

\begin{table}
\begin{center}
\begin{tabular}{||l||c|c|c|||l|c|c|c||} \hline \hline
Knot & Tor & Rank & Col type & Knot & Tor  & Rank & Col type
\\ \hline \hline
$3_1$ & 3  &  3 & 3 Trivial   &$9_{11}$   &33 &3 & 3 Trivial  \\ \hline
            & 0 &4 & 6 Non-trivial &   & 0 & 4 & 6 Non-trivial   \\ \hline
$6_1$ & 9  & 3 & 3 Trivial  &$9_{15} $ & 39 & 3 & 3 Trivial    \\ \hline
            & 0 & 4 & 6 Non-trivial  &          & 0 & 4 & 6 Non-trivial  \\
\hline
$7_4$  & 15 & 3 & 3 Trivial    &$9_{16} $  & 39 &3&3 Trivial   \\ \hline
             &  0 & 4 &6 Non-trivial & & 2& 4&6 Non-trivial    \\ \hline
$7_7$    & 21 & 3 & 3 Trivial& $9_{17} $ & 39 & 3 &3 Trivial \\ \hline
                & 0  & 4 & 6 Non-trivial &    & 0& 4  &6 Non-trivial  \\ \hline
$8_5$  & 21& 3& 3 Trivial & $9_{23}$  & 45 & 3 &3 Trivial \\ \hline
             & 2 & 4 & 6 Non-Trivial  & & 0 & 4 &6 Non-trivial \\ \hline
$8_{10}$ & 27 &3 & 3 Trivial  &$9_{24}$ & 45 & 3 & 3 Trivial  \\ \hline
                   & 2 & 4&  6 Non-trivial  &   & 2 & 4&6 Non-trivial  \\
\hline
$8_{11}$   &27 &3 &3 Trivial   &$9_{28}$ &51 &3 &3 Trivial \\ \hline
& 0 & 4 & 6 Non-trivial &     & 2 & 4 &6  Non-trivial  \\ \hline
$8_{15}$   &33 &3 &3 Trivial   & $9_{29}$ &51 &3 &3 Trivial \\ \hline
                    & 2& 4&6 Non-trivial&    & 0 & 4 &6 Non-trivial   \\ \hline
$8_{18}$    & 3, $\;$ 15 &3 & 3 Trivial& $9_{34}$ & 69 &3  & 3 Trivial \\
\hline
                     & 3 & 4 & 24 Non-trivial  &   & 0 & 4 & 6 Non-trivial  \\
\hline
$8_{19}$    &3 &3 & 3 Trivial    &$9_{35}$ & 3, $\;$ 9 &3 & 3 Trivial \\ \hline
                     & 2 & 4 &6  Non-Trivial &             & 0 & 4 & 18
Non-trivial   \\ \hline
$8_{20} $   & 9&3&3 Trivial    &      & 2 & 4 & 6 Non-trivial \\ \hline
            & 2&4 & 6 Non-trivial & $9_{37}$ & 3, $\;$ 15 & 3 & 3 Trivial   \\
\hline
$8_{21}$    &15&3 &3 Trivial  & & 2& 4 & 6 Non-trivial \\ \hline
                     &2 & 4 & 6 Non-trivial &   &   0& 4 & 18 Non-trivial \\
\hline
$9_1$   &9 &3 &3 Trivial & $9_{38}$ & 57&3& 3 Trivial  \\ \hline
              & 0 & 4 & 6 Non-trivial   & & 0 & 4 & 6 Non-trivial \\ \hline
$9_{2}$& 15 &3&3 Trivial  & $9_{40}$ & 5, $\;$ 15 & 3 & 3 Trivial \\ \hline
              & 0 & 4 & 6 Non-trivial &  & 4 & 4 & 6 Non-trivial \\ \hline
 $9_4$  &  21 &3& 3  Trivial   & $9_{46}$ & 3, $\;$ 3 & 3 & 3 Trivial     \\
\hline
             & 0 & 4& 6 Non-trivial   & & 0 & 4 & 18 Non-trivial \\ \hline
  $9_6$   & 27 & 3  & 3 Trivial    & & 2 & 4 & 6 Non-trivial \\ \hline
         & 0 & 4 & 6 Non-trivial   &$9_{47}$ & 3, $\;$ 9 & 3 & 3 Trivial \\
\hline
$9_{10}$ &  33 & 3 & 3 Trivial & & 0 &4 & 24 Non-trivial     \\ \hline
                   & 0 & 4 & 6 Non-trivial  & $9_{48}$ & 3, $\;$ 9 & 3 &3
Trivial \\ \hline
& & &    &   & 0 & 4 & 18 Non-trivial  \\ \hline
& & &&      & 2 & 4 & 6 Non-trivial \\ \hline
\end{tabular}\end{center}
\caption{A table of module invariants for $3$-colorable knots}
\label{qmoduletable3}
\end{table}

\begin{table}[h]
\begin{center}
{\begin{tabular}{||l||c|c|c|||l|c|c|c||} \hline \hline
Knot & Tor & Rank & Col type & Knot & Tor & Rank & Col type
\\ \hline \hline
$4_1$ & 5, $\;$ 5 &5 & 5 Trivial &
 $9_2$& 15, $\;$ 15 & 5 & 5 Trivial \\ \hline
             &0 & 7 & 20 Non-trivial &     
       & 0 & 7 & 20 Non-trivial \\ \hline
$5_1$ & 5, $\;$ 5 & 5 & 5 Trivial &
 $9_{12}$& 35, $\;$ 35 & 5 & 5 Trivial \\ \hline
              & 0 & 7 &20 Non-trivial &      
       & 0 & 7 & 20 Non-trivial \\ \hline
$7_4$ & 15, $\;$ 15 & 5& 5 Trivial & 
$9_{23}$& 45, $\;$ 45 & 5 & 5 Trivial \\ \hline
             & 0 & 7 &20 Non-trivial &     
        & 0 & 7 & 20 Non-trivial \\ \hline
$8_7$ & 25, $\;$ 25 & 5 & 5 Non-trivial &
 $9_{24}$& 45, $\;$ 45 & 5 & 5 Trivial \\ \hline
             & 0 & 7 & 20 Trivial&          
   & 0 & 7 & 20 Non-trivial \\ \hline
$8_8$& 25, $\;$ 25 & 5 & 5 Trivial &
$9_{31}$& 55, $\;$ 55 & 5 & 5 Trivial \\ \hline
             & 0 & 7 & 20 Non-trivial &       
      & 0 & 7 & 20 Non-trivial \\ \hline
$8_{16}$& 35, $\;$ 35 & 5 & 5 Trivial & 
$9_{37}$& 3, $\;$ 3, $\;$ 15, $\;$ 15 & 5 & 5 Trivial \\ \hline
             & 0 & 7 & 20 Non-trivial &   
          & 0 & 7 & 20 Non-Trivial \\ \hline
$8_{18}$& 3, $\;$ 3, $\;$ 15, $\;$ 15 & 5 & 5 Trivial &
$9_{39}$ & 55, $\;$ 55 & 5 & 5 Trivial \\ \hline
             &  2, $\;$ 2 &7 & 20 Non-trivial & 
            & 0 & 7 & 20 Non-trivial \\ \hline
$8_{21}$& 15, $\;$ 15 & 5 & 5 Trivial&
 $9_{40}$& 5, $\;$ 5, $\;$ 15, $\;$ 15 & 5 & 5 Trivial \\ \hline
             & 0 & 7 & 20 Non-trivial&   
          & 5 & 7 & 120 Non-Trivial \\ \hline
 & & & &  $9_{48}$& 5, $\;$ 5, $\;$ 5, $\;$ 5 & 5 & 5 Trivial \\ \hline
 & & & &      & 0 & 7 & 120 Non-trivial \\ \hline
\end{tabular} } \end{center}
\caption{A table of module invariants for $5$-colorable knots}
\label{qmoduletable5}
\end{table}


\begin{remark} {\rm
The construction of the quandle module invariant 
is similar to the construction of Alexander modules from Burau 
representation. Also, the cokernel appears in the 
definition of the Bowen-Franks~\cite{BF}
groups for symbolic dynamical systems. Thus  relations to 
covering spaces, as well as  dynamical systems related to braid groups, 
are expected.

In \cite{Lin,Wada} group representations of knot groups are used to
define twisted Alexander polynomials. In that situation, the
representation can be viewed as a matrix (which depends on the image of
the meridian) assigned to a crossing. The quandle module invariant is
related when the quandle colorings are given by knot group
representations. 
The general relation can be understood via Fox calculus.

The quandle module
invariants might give lower bounds to the number of strands needed in a
braid representation of a knot.
Future studies will include more detailed investigations of this
invariant, the non-abelian cocycle invariant (Section 5), and the
generalized cocycle invariant of classical knots (Section 6). Our main
purpose in the current paper is to indicate the strength of the
generalized cocycle invariants  for knotted surfaces (Section 7).
} \end{remark}

\section{Knot invariants  from non-abelian $2$-cocycles}
\label{nonabsec}

\subsection*{Non-abelian $2$-cocycles}
Let $X$ be a quandle and $H$ a (not necessarily abelian) group. 
A function $\beta: X  \times X \rightarrow H$ 
is a {\it rack $2$-cocycle} \cite{AG}  if
$$ \beta(x_1, x_2) \beta(x_1*x_2, x_3) 
= \beta(x_1, x_3) \beta (x_1 * x_3 , x_2 * x_3 ) $$
is satisfied for any $x_1$, $x_2$, $x_3 \in X$.
If a rack $2$-cocycle further satisfies  $\beta(x,x)=1$ 
for any $x \in X$, then it is called 
a {\it quandle $2$-cocycle}  \cite{AG}. 
The set of quandle  $2$-cocycles is denoted by $Z^2_{\rm CQ}(X;H)$.
Two cocycles $\beta$, $\beta'$ are {\it cohomologous} if 
there is a function $\gamma: X \rightarrow H $ such that 
$$ \beta'(x_1, x_2) = \gamma(x_1)^{-1}\beta(x_1, x_2) \gamma(x_1 * x_2) $$
for any $x_1, x_2 \in X$. 
An equivalence class is called a {\it cohomology class}.
The set of cohomology classes is denoted by $H^2_{\rm CQ}(X;H)$.
These definitions agree with those in \cite{CJKLS} if $H$ is
an abelian group. When $H$ is not necessarily abelian,
the $2$-cocycles $\beta$ are called (constant)
 {\it  $2$-cocycles} in 
\cite{AG}.  We call such a $2$-cocycle 
{\it non-abelian } when $H$ is not an abelian group.

Let $S$ be a set and $\Sym_S$ denotes the permutation group on $S$.
Let $E=S \times X$ and $\beta \in Z^2_{\rm CQ}
(X; \Sym_S)$.
Then the binary operation on $E$
$$ (a_1, x_1) * (a_2 , x_2)=( a_1 \cdot \beta(x_1, x_2) , x_1 * x_2)$$
defines a quandle structure on $E$. We call this $E$ the non-abelian 
extension of $X$ by $\beta$, and denote it by
$E=E(X,S, \beta)$.

A quandle $X$ is {\it decomposable} \cite{AG} if it is a disjoint union 
$X=Y \sqcup Z$ such that $Y*X=Y$ and $Z*X=Z$,
where $Y*X=\{ y*x\  | \ y \in Y, \ x \in X \} $.
A  quandle is {\it indecomposable} if it is not decomposable.
For a rack $X$, let $\phi_y : X \rightarrow X$ be the quandle 
isomorphism defined by $\phi_y (x)=x*y$, $x \in X$. 
The subgroup $\mbox{Inn}(X)$ of the group $\mbox{Aut}(X)$ of 
automorphisms of $X$ generated by $\phi_y$, $y \in X$, is 
called the inner automorphism group. The same groups are 
defined for quandle automorphisms for a quandle $X$.

\begin{lemma} \label{AGconstlemma} {\bf \cite{AG}} 
Suppose  $Y$ is indecomposable and 
$f: Y \rightarrow X$ is a quandle surjective homomorphism
such that 
$\phi_x = \phi_y$  
if  $f(x)=f(y)$ for $x, y \in Y$.
Then $Y$ is a non-abelian extension $S \times_{\beta} X$ 
for some $\beta \in Z^2_{\rm CQ}(X;\Sym_S )$.
\end{lemma}

The following construction was given in \cite{AG} (Example 2.13).
Let $X$ be an indecomposable finite rack, $x_0 \in X$ a fixed element,
$G=\mbox{Inn}(X)$, and $H=G_{x_0}$ be the subgroup of $G$ 
whose elements fix $x_0$. Let $S$ be a finite set and 
$\rho : H \rightarrow \Sym_S$ be a group homomorphism. There is a bijection 
$G/H \rightarrow X$ given by 
$g \mapsto g(x_0)$. 
Fix a set-theoretic section $s: X \rightarrow G$. 
Thus $s(x) \cdot x_0=x$ for all $x \in X$. 

\begin{lemma} \label{AGcocylemma} {\bf \cite{AG}} 
The element $t(x,y)=s(x) \phi_{y} s(x*y)^{-1} $ is in $H$
for any $x, y \in X$, 
and $\beta(x,y)=\rho (t(x,y))$ is a rack $2$-cocycle. 

The $2$-cocycle $\beta$ is a quandle $2$-cocycle if and only if
$\rho(\phi_{x_0})=1 \in \Sym_S$.
\end{lemma}

\subsection*{Definitions}

We  define a new cocycle invariant using the  non-abelian cocycles.
Let $L=K_1 \cup \cdots  \cup K_r$ be a classical oriented
link diagram on the plane,
where $K_1, \ldots, K_r$ are connected components, 
for some positive integer $r$.
Let ${\cal T}_i$, for, $i=1, \ldots, r$, be  the set of crossings such that
the under-arc is from the component $i$.

Let $X$ be a quandle, $H$ a group, $\beta \in Z^2_{\rm CQ}(X;H)$.
Let ${\cal C} \in \mbox{Col}_X(L)$ be a coloring of $L$ by $X$.
Let $( b_1, \ldots , b_r)$ be the set of base points on the 
components $(K_1, \ldots, K_r)$, respectively.
Let $(\tau^{(j)}_1, \ldots, \tau^{(j)}_{k(j)} )$ be the crossings in
${\cal T}_j$, $j=1, \ldots, r$, that appear in this order 
when one starts from $b_j$ and travels $K_j$ in the given orientation.

At a crossing $\tau$, let $x_{\tau}$ be the color on the under-arc 
from which the normal of the over-arc points; let $y_{\tau}$ be 
the color on the over-arc. The {\it Boltzmann weight} at $\tau$
is $B(\tau, {\cal C})= \beta(x_{\tau}, y_{\tau})^{\epsilon(\tau)}$, where 
$\epsilon(\tau)$ is $\pm 1$ depending on 
whether $\tau$ is positive or negative,
respectively. 
For a group element $h \in H$, denote by $[h]$ the conjugacy class 
to which $h$ belongs.

\begin{definition} \label{cocydef} {\rm 
The family of vectors of conjugacy classes
$$
 \vec{\Psi} (L)= \vec{\Psi}_{(X, H, \beta)} (L) 
= ( [ \Psi_1 (L, {\cal C}) ] , \ldots ,  [ \Psi_r (L, {\cal C}) ] )_{  {\cal C } \in {\rm Col}_X(L) } 
$$ where 
$$ \Psi_i (L, {\cal C})=   \prod_{j=1}^{k(i)} \beta(x_{\tau^{(i)}_{j}} , 
y_{\tau^{(i)}_{j}})^{\epsilon(\tau^{(i)}_j )}  ,  $$
is called the {\it  conjugacy quandle cocycle invariant } of a link.
} \end{definition}

These  cocycle invariants include abelian cocycle
invariants defined in \cite{CJKLS} as a special case 
(when $H$ is abelian).

\begin{figure}[h]
\begin{center}
\mbox{
\epsfxsize=2in
\epsfbox{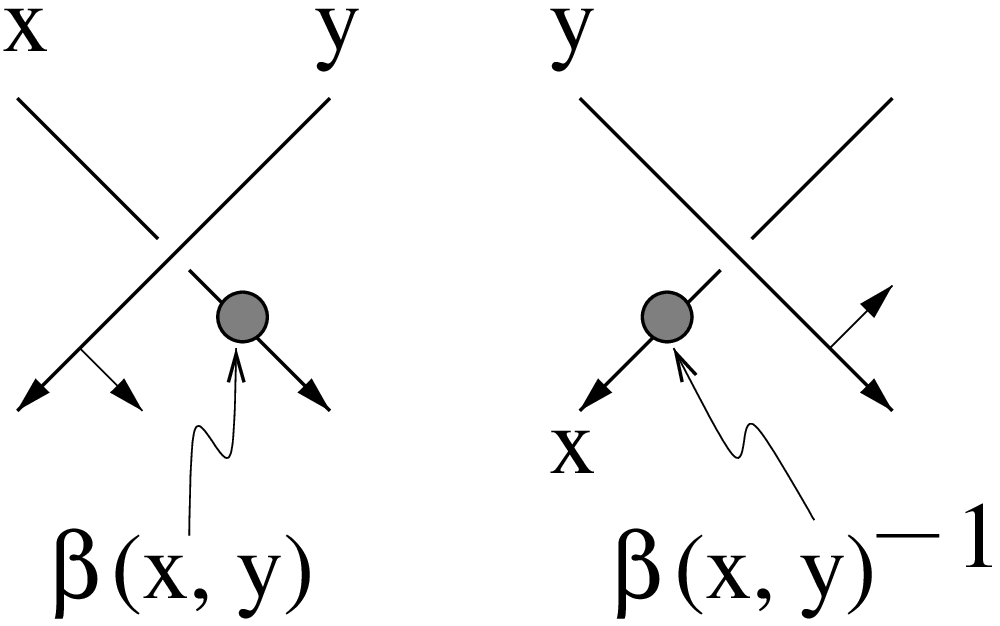} 
}
\end{center}
\caption{ The beads interpretation of the invariant  }
\label{beads} 
\end{figure}

\begin{figure}[h]
\begin{center}
\mbox{
\epsfxsize=3in
\epsfbox{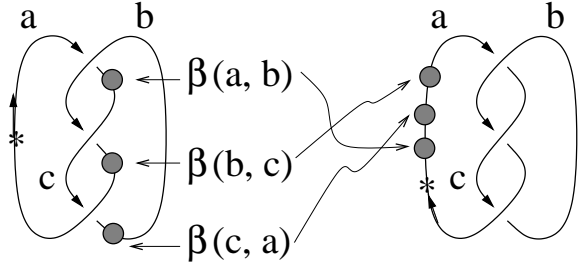} 
}
\end{center}
\caption{ Making a beads necklace with trefoil }
\label{beadstre} 
\end{figure}

\begin{remark} {\rm
The invariant has the following interpretation of sliding beads 
along knots, similar to Hopf algebra invariants defined in \cite{KauRad}.
Let a knot diagram $K$ and its coloring ${\cal C}$ by 
a finite quandle $X$ be given, and let $\Phi(K)$ be the cocycle
invariant with a  $2$-cocycle $\beta \in Z^2_{\rm CQ}(X;H)$ 
for a (non-abelian) coefficient group $H$. 

Put a bead on the underarc just below each crossing as shown
in Fig.~\ref{beads}. Each bead is assigned 
the weight at the crossing, as in the left of the figure
for a positive crossing 
and as in the right for a negative crossing, respectively.
This process is depicted in Fig.~\ref{beadstre} in the left
for the trefoil.

Pick a base point (which is depicted by $*$ in the figure) on the 
diagram, and push it in the given orientation of the knot. 
With it the beads are pushed along in the order the base point
encounters them. When the base point comes back to near
the original position, it has collected all the beads.
This situation is depicted in Fig.~\ref{beadstre} in the right.
The beads, read from the base point, are aligned 
in the order $\beta(a,b)$, $\beta(b,c)$, and $\beta(c,a)$ 
in this order in the figure, and 
the conjugacy class 
$[\beta(a,b) \beta(c,a) \beta(b,c)]$  is the contribution to the invariant 
for this coloring. 
} \end{remark}

\begin{theorem} 
The quandle cocycle invariant 
$\vec{\Psi}$ is well defined. 
Specifically, let $L_1, L_2$ be two link diagrams of ambient isotopic links,
and $\vec{\Psi}(L_1), \vec{\Psi}(L_2)$ be their quandle cocycle invariants.
Then there is a bijection 
$\eta:  \vec{\Psi}(L_1) \rightarrow \vec{\Psi}(L_2)$. 
\end{theorem}
{\it Proof.\/} 
The fact that $\vec{\Psi}$ does not change 
by Reidemeister moves for fixed base points is similar
to the proof in \cite{CJKLS} for the knot case and that in \cite{CENS}
for the link case, proved for abelian cocycle invariants, 
except one just observes that the order of group elements is also preserved 
under Reidemeister moves.

\begin{figure}[h]
\begin{center}
\mbox{
\epsfxsize=2.5in
\epsfbox{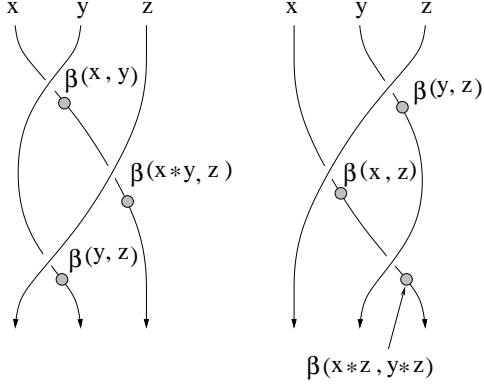} 
}
\end{center}
\caption{ The beads and the type III move  }
\label{beadstype3} 
\end{figure}

Specifically, the changes under the type III Reidemeister move is 
depicted in Fig.~\ref{beadstype3}. 
On the bottom string (that goes from top left to bottom right),
two cocycles are assigned in both of the left and right
of the figure, and they are equal with the given order
(from top to bottom), by the $2$-cocycle condition. 
The only cocycle assigned on the middle string has the 
same value for the left and right of the figure, 
and they also occupy the same position when the cocycles are 
read (as beads are slidden) along the component.
Thus the ordered elements do not change by this move.

A change of base points causes cyclic permutations
of Boltzmann weights, and hence the invariant is defined up to conjugacy.
\qed

\begin{figure}[h]
\begin{center}
\mbox{
\epsfxsize=2in
\epsfbox{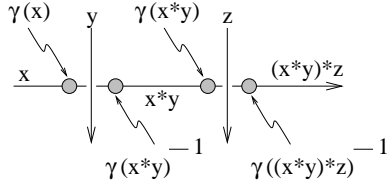} 
}
\end{center}
\caption{Canceling beads terms  }
\label{beadscancel} 
\end{figure}

\begin{proposition} \label{cobprop}
The  cocycle invariant $\vec{\Psi}$ is trivial 
(i.e., $\vec{\Psi}$ consists of vectors whose entries are the conjugacy class
of the identity element)  if the 
cocycle used is a coboundary.
\end{proposition}
{\it Proof.\/} 
The proof is similar to the proof of an analogous theorem 
 in \cite{CJKLS}.
If $\beta \in Z^2_{CQ}(X;H)$ is null-homologous, 
then $\beta(x,y)=\gamma(x) \gamma(x*y)^{-1}$
for some function $\gamma: X \rightarrow H$.
Let $\beta(x,y)$, $\beta(x*y, z)$ be  consecutive $2$-cocycles 
in $\Psi_i(K, {\cal C})$, a contribution to the cocycle
invariant for a coloring ${\cal C}$ for the $i\/$th component.
Then in  $\Psi_i(K, {\cal C})$, they form a product 
$\cdots \beta(x,y) \beta(x*y, z) \cdots$, which is equal to
$\cdots (\gamma(x) \gamma(x*y)^{-1}) (  \gamma(x*y) \gamma( (x*y)*z )^{-1} )\cdots$,
and the middle terms cancel. The left and the right terms 
cancel with the next adjacent  terms, and obtain $\Psi_i(K, {\cal C})=1$.
This proves the proposition. The situation is 
easier to visualize diagrammatically. The bead representing 
$\beta(x,y)$ is represented by two separate ordered beads representing 
$\gamma(x)$ and $\gamma(x*y)^{-1}$ as depicted at the left 
crossing in Fig.~\ref{beadscancel}. 
Thus all the beads cancel after going around each component once.
\qed

Recall that an action of a group H on a set X is said \emph{free}
if for all $h\in H$, $x\in X$ one has that $h\cdot x=x$ implies $h=1$.

\begin{proposition} \label{colorextprop}
Let $H \subseteq \Sym_S$ and $\beta \in  Z^2_{\rm CQ}(X;H)$.
If the invariant $\vec{\Psi}(L)$ is trivial, then
any coloring of $L$ by $X$ extends to a coloring by $E(X,S,\beta)$.
Conversely, if any coloring of $L$ by $X$ extends to a coloring
by $E(X,S,\beta)$, then  $\vec{\Psi}(L)$ is trivial,
provided that the action on $S$ of the subgroup
of $H$ generated by the image of $\beta$ is free.
\end{proposition}
{\it Proof.\/}
This  proof is similar to the proof of the
corresponding theorem in \cite{CENS}.
Pick an element $s_0 \in S$ and color the arc with the base
point $b_i$ on the $i\/$th component by $(s_0, x) \in E=E(X, S, \beta)$,
where $E$ is identified with $S \times X$.
Go along the component, and after passing under the first crossing,
the color of the second arc is required to be
$(s_0 \beta(x, y), x*y) $, if the over-arc is colored by $y \in X$.
Continuing this process, when we come back to the base point,
the color is $(s_0 \Psi_i(L, {\cal C}), x)$,
and the proposition follows.
\qed

The freeness requirement is satisfied, for example, if
$S=H$ and $H$ acts as multiplication on the right.

When a $2$-cocycle constructed by the method of Lemma~\ref{AGcocylemma}
is used to define the invariant, we have the following interpretation
that is an analogue of Propositions~\ref{cobprop} and \ref{colorextprop}. 
Let $L=K_1 \cup \cdots \cup K_r$ be an oriented link diagram,
and $b_1, \ldots, b_r$ arbitrarily chosen and fixed 
 base points on each component.
Let $\beta$ be a $2$-cocycle constructed in Lemma~\ref{AGcocylemma}, 
so that $\beta(x,y)=\rho (t(x,y)) $ for any $x, y$ in a quandle $X$,
and $t(x,y)=s(y) \phi_{y} s(x*y)^{-1}$.
{}From the proof of Proposition~\ref{colorextprop},
the contribution to the cocycle invariant is written as 
\begin{eqnarray*} 
\Psi_i(L, {\cal C}) &= & 
 \beta(x_{\tau^{(i)}_{1}} , 
y_{\tau^{(i)}_{1}})^{\epsilon(\tau^{(i)}_1 )} 
 \cdot
 \beta(x_{\tau^{(i)}_{2}} , 
y_{\tau^{(i)}_{2}})^{\epsilon(\tau^{(i)}_2 )}
\cdots 
 \beta(x_{\tau^{(i)}_{k(i)}} , 
y_{\tau^{(i)}_{{k(i)}}})^{\epsilon(\tau^{(i)}_{k(i)} )}  \\
&=& (\rho ( s(x_{\tau^{(i)}_{1}}) \phi_{y_{\tau^{(i)}_{1}}}
 s(x_{\tau^{(i)}_{2}} )^{-1}) ) 
\cdot
(\rho ( s(x_{\tau^{(i)}_{2}} ) \phi_{y_{\tau^{(i)}_{2}}}
 s(x_{\tau^{(i)}_{3}} )^{-1}) ) 
\cdots
(\rho ( s(x_{\tau^{(i)}_{k(i)}}) \phi_{y_{\tau^{(i)}_{k(i)}}}
 s(x_{\tau^{(i)}_{k(i)}} )^{-1}) ) \\
&=& \rho (  s(x_{\tau^{(i)}_{1}})
\cdot (\phi_{y_{\tau^{(i)}_{1}}}  \phi_{y_{\tau^{(i)}_{2}}}
\cdots  \phi_{y_{\tau^{(i)}_{k(i)}}} )
\cdot  s(x_{\tau^{(i)}_{k(i)}})^{-1} )
\end{eqnarray*}
Thus we obtain the following.

\begin{lemma} \label{longilemma}
If a $2$-cocycle constructed in Lemma~\ref{AGcocylemma}
is used to define the conjugacy cocycle invariant, 
then the contribution to the $i\/$th component 
$\Psi_i(L, {\cal C})$ is computed by 
$$\Psi_i(L, {\cal C}) = \rho (  s(x_{\tau^{(i)}_{1}})
\cdot (\phi_{y_{\tau^{(i)}_{1}}}  \phi_{y_{\tau^{(i)}_{2}}}
\cdots  \phi_{y_{\tau^{(i)}_{k(i)}}} )
\cdot  s(x_{\tau^{(i)}_{1}} )^{-1} ) . $$
\end{lemma}

The expression 
$y_{\tau^{(i)}_{1}}\cdots y_{\tau^{(i)}_{k(i)}}$
is the  sequence  of 
colors that one encounters traveling along the component,
picking up $y_{\tau^{(i)}_{j}}$ as one goes under the $j\/$th 
crossing. Thus this sequence corresponds to the longitudinal
element in the fundamental group.

\subsection*{Constructions}

We follow Lemma~\ref{AGcocylemma} to construct explicit examples
of non-abelian cocycle invariants from conjugacy classes of groups.
Let $X$ be a conjugacy class in a finite group and let $G=\langle X \rangle$
be the subgroup generated by $X$. 
Let $x_0 \in G$ be a fixed element and 
$H=Z_{x_0} =\{ x \in G \ | \ x_0 x = x x_0 \} $.
Let $\alpha: G \rightarrow \mbox{Inn}(X)$ be 
the map induced from the conjugation 
$g \mapsto ( x \mapsto gxg^{-1})$. 
Then the kernel of $\alpha $ is the center $Z(G)$ 
and $\mbox{Inn}(X)$ is isomorphic to $G / Z(G)$. 

To evaluate cocycle invariants for specific examples, we take 
 $X=(2,1,\ldots, 1)$ ($(n-2)$ copies of $1$),
the conjugacy class consisting of 
transpositions of $\Sym_n$, $n \geq 5$, then 
$G=\langle X \rangle=\Sym_n$, and $\mbox{Inn}(X)=G$.
Let $x_0=(1\ 2)$, then 
$$H=\langle (1 \ 2), (i \ j ) \ | \ 2 < i < j \leq n \rangle \cong 
\Z_2 \times \Sym_{n-2}, $$
and  $N(x_0)$, the normal closure of $x_0$ in $H$, 
is equal to $\langle (1 \ 2) \rangle=\{ 1, (1\ 2)\}$, 
so that  $H/N(x_0) \cong \Sym_{n-2}$, where 
$\Sym_{n-2}$ is identified with the symmetric group on 
$n-2$ letters $\{ 3, 4, \ldots, n \}$.
Thus take $S=\{  3, 4, \ldots, n \}$ and regard 
$H/N(x_0)$ as  $\Sym_S$.
Let the map 
$\rho : H \rightarrow H / N(x_0) = \Sym_S$ be the projection.
Then $\rho(x_0)=1$, so that the condition for a quandle 
cocycle is satisfied. 
Let $s:  X \rightarrow G$ be defined by $s(i\ j)=(1\ i)(2 \ j)$,
then $s$ defines a section.
By Lemma~\ref{AGcocylemma}, this set up defines 
a $2$-cocycle $\beta \in Z^2_{\rm CQ}(X; \Sym_S)$. 

\begin{example} \label{cocyex}
{\rm
Let $n=5$, then, for example, $\beta( (1\ 4) , (2 \ 3) )$ can be computed 
as 
\begin{eqnarray*} 
\beta( (1\ 4) , (2 \ 3) ) & = & \rho ( t( (1\ 4) , (2 \ 3) ) ) \\
&=& \rho( s( (1\ 4)) \phi_{   (2 \ 3) } s(  (1\ 4)* (2 \ 3) )^{-1}) \\
&=& \rho( (2\ 4)  (2 \ 3) [ s( (2 \ 3) (1\ 4)(2 \ 3)) ]^{-1}) \\
&=& \rho(  (2\ 4)  (2 \ 3) (2\ 4) ) \\
&=& \rho( (3\ 4) ) \\
&=& (3\ 4) .
\end{eqnarray*}

} \end{example}

\begin{example} {\rm 
We evaluate the cocycle invariant using the preceding $2$-cocycle
for a Hopf link $L=K_1 \cup K_2$ depicted in Fig.~\ref{hopf}.

\begin{figure}[h]
\begin{center}
\mbox{
\epsfxsize=2in
\epsfbox{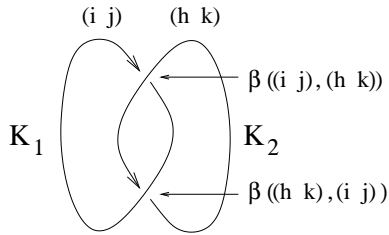} 
}
\end{center}
\caption{Hopf link  }
\label{hopf} 
\end{figure}

Let $n=5$ as in Example~\ref{cocyex}.
There are two arcs in a Hopf link, also denoted by $K_1$ and $K_2$.
In the figure, the colors $(i\ j)$ and $(h\ k)$
are assigned. Consider the case $(i\ j)=(1\ 4)$ and $(h\ k)=(2\ 3)$,
which certainly defines a coloring ${\cal C}$ of $L$. 
By the computation in Example~\ref{cocyex}, 
the contribution to the invariant $\vec{\Psi}(L)$ is
$$\Psi_1(L, {\cal C}) =
\beta( (1\ 4) , (2 \ 3) )=(3\ 4)=\Psi_2(L, {\cal C}), $$
which contributes
the pair $([ (3\ 4)], [(3\ 4)])$ of conjugacy classes of 
$\Sym_3=\Sym_{ \{ 3,4,5 \} }$.

} \end{example}

\subsection*{The non-abelian cocycle invariant is a quantum invariant}

Let $X$ be a finite quandle, $G$ a finite group and
$\beta:X\times X\to G$ a $2$-cocycle. 
Let $\vec{\Psi}(L)=([\Psi_1(L, {\cal C})], \ldots, [\Psi_r(L, {\cal C})])$
be the above defined invariant. 

Let $\rho:G\to\Aut(V)$ be a representation of $G$.
Let $W=\CC X\otimes V$ with the braiding given by
$c((x\otimes v)\otimes(y\otimes w))
	= (y\otimes w)\otimes(x*y\otimes v\cdot\beta(x,y))$.
Since $X$ and $G$ are finite, $W$ is a (right--right) Yetter-Drinfeld
module over some finite group,
(see \cite{AG}) 
and we can consider the quantum link invariants
coloring links by $W$ (the ribbon structure is the identity map,
 see \cite{Grana}). Denote this invariant by  $\Psi(L,W)$.
Then we have the following, an analogue of \cite{Grana}.

\begin{proposition}
$\tr(\rho\vec{{\Psi}}(L) )  = \Psi(L,W).$ 
(Here $\tr(\rho\vec{{\Psi}}(L))=\prod_{i=1}^r\tr([\Psi_i(L, {\cal C})])$
is the trace in the tensor product via the inclusion $\Aut(V)\times\cdots\times\Aut(V)%
\hookrightarrow\End(V\otimes\cdots\otimes V)$ ($r$ components)). 
\end{proposition}
{\it Proof.\/} 
Consider $L$ as the closure of a braid
$z=\sigma_{i_1}\cdots\sigma_{i_l}\in\BB_n$.
Let $\bar z\in\Sim_n$ be the projection of $z$ in the symmetric group.
To compute $\Psi(L,W)$ we take the trace of the map defined by
$$z:W^{\epsilon_1}\otimes\cdots\otimes W^{\epsilon_n}
	\to W^{\epsilon_{\bar z(1)}}\otimes\cdots\otimes W^{\epsilon_{\bar z(n)}},$$
where $\epsilon_i=\pm 1$ and $W^{-1}$ stands for $W^*$. To compute
this trace we take the basis $X$ of $\CC X$, a basis $\{v_i\ |\ i\in I\}$
of $V$ and the basis $\{x\otimes v_i\ |\ x\in X, i\in I\}$. Notice
that the first components of the braiding in $W$ are given by the quandle $X$.
Thus, if we apply the map $z$ to an element of the form
$(x_1\otimes v_{i_1})\otimes\cdots\otimes(x_n\otimes v_{i_n})$, this element
will not contribute to the invariant unless we receive $x_1,\ldots,x_n$
as the first components; i.e. if $x_1,\ldots,x_n$ defines a coloring of $L$.
Consider therefore the colorings of $L$. Take a particular coloring given
by $x_1,\ldots,x_n$ at the bottom of the braid. For any $n$-tuple of elements
$v_{i_1},\ldots,v_{i_n}$ in the basis of $V$, we apply successively
$\sigma_{i_1},\ldots,\sigma_{i_l}$ and notice that the $v_i$'s change by
$\beta$ exactly in the same way as defined in the invariant 
$\overline{{\Psi}}$. 
Now, 
$\tr(\rho\vec{{\Psi}})$ 
is the sum over the colorings of traces of
products of $\rho(\beta_{x,y})$ for several $x,y$'s, and $\Psi$ is
the sum over colorings of traces of maps made by $\beta_{x,y}$ for
the same pairs $x,y$ as for 
$\tr(\rho\vec{{\Psi}})$. 
The last step is to
notice that the order of the $\beta$'s used in the definition of 
$\vec{{\Psi}}$ 
is the same one that one gets by taking traces, thanks to the following remark:
Let $f_1,\ldots,f_n\in\End(V)$ be linear maps and consider the map
$F\in\End(V^{\epsilon_1}\otimes\cdots\otimes V^{\epsilon_n})$ given by
$$F(x_1^{\epsilon_1}\otimes\cdots\otimes x_n^{\epsilon_n})
	= \bar z(f^{\epsilon_1}(x_1^{\epsilon_1})\otimes\cdots\otimes
		f^{\epsilon_n}(x_n^{\epsilon_n})).$$
Write $\bar z=(i_1,\ldots,i_{a_1})(i_{a_1+1},\ldots,i_{a_1+a_2})\cdots
	(i_{a_1+\cdots+a_{r-1}+1},\ldots,i_{a_1+\cdots+a_r})$ the decomposition of
$\bar z$ in disjoint cycles.
Then
$$\tr F=\prod_{j=1}^r\tr(f_{a_1+\cdots+a_{j-1}+1}^{\epsilon_{a_1+\cdots+a_{j-1}+1}}
	f_{a_1+\cdots+a_{j-1}+2}^{\epsilon_{a_1+\cdots+a_{j-1}+2}}
	\cdots f_{a_1+\cdots+a_{j}}^{\epsilon_{a_1+\cdots+a_{j}}}).$$
This is easy to see, though the notation is cumbersome. Let us illustrate it
with the easy example $z=\sigma_1\in B_2$; we must prove then that
$$\tr((v\otimes w)\mapsto (g(w)\otimes f(v)) = \tr (fg),$$
but the left hand side of the equation is $\sum_{i,j}g_{j,i}f_{i,j}$
(here $f_{i,j}$ and $g_{i,j}$ are the matrix coefficients of $f$ and $g$
in the basis $\{v_i\}$ of $V$) and this coincides with the right hand side.
\qed

\section{Knot invariants  from generalized $2$-cocycles} \label{cocyinvsec}

\subsection*{Definitions}

Let $X$ be a finite quandle and $\Z (X)$ be its quandle algebra
with generators $\{ \eta _{x,y}^{\pm 1}  \}_{ x, y \in X} $
and $\{ \tau _{x,y}  \}_{ x, y \in X} $.
Recall (Example~\ref{AGexample}) that the  enveloping group
 $G_X$ of a quandle $X$ is a quotient 
of the free group $F(X)$ on $X$ defined by
 $G_X=\langle x \in X \ | \ x*y=yxy^{-1} \rangle$.
Then 
  $\Z (X)$ maps onto the group algebra $\Z G_X$ 
by   $\eta_{x,y } \mapsto y$, $\tau_{x,y} \mapsto 1-x*y$
so that any left $\Z G_X$-module has a structure of a $\Z (X)$-module
\cite{AG}. Let $G$ be an abelian group with 
this  $\Z (X)$-module structure.
Let $\kappa_{x,y}$ be a generalized quandle $2$-cocycle of $X$ with 
the coefficient group $G$. 
Thus the generalized $2$-cocycle condition, in this setting, is 
written as 
$$ z \kappa_{x,y} + \kappa_{x*y, z} + ((x*y)*z) \kappa_{y,z}
=  \kappa_{y,z} + (y*z)  \kappa_{x, z} + \kappa_{x*z,y*z} .$$
We define a cocycle invariant using this $2$-cocycle. 

A knot diagram $K$ is given on the plane. 
Recall that it is oriented, and has orientation normals. 
There are four regions near each crossing, divided by the arcs of the diagram.
The unique 
region into which 
both normals (to over- and under-arcs) point
is called the {\it target} region.  
Let $\gamma$ be an arc from the region at infinity of the plane
 to the target region of a given crossing $r$, that intersect 
$K$ in  finitely many points transversely, missing crossing points. 
Let $a_i$, $i=1, \ldots, k$, in this order, be the arcs of $K$
that intersect $\gamma$ from the region at infinity to the crossing.
Let ${\cal C}$ be a coloring of $K$ by a fixed finite quandle $X$.

\begin{definition}
{\rm The {\em Boltzmann weight} $B( {\cal C}, r, \gamma)$
for the crossing $r$, for a coloring ${\cal C}$, with respect to $\gamma$,
is defined by 
$$ B( {\cal C},r, \gamma)=
 \pm  ({\cal C}(a_1)^{\epsilon (a_1) } {\cal C}( a_2)^{\epsilon (a_2) }
 \cdots{\cal C}( a_k)^{\epsilon (a_k) } ) 
 \kappa_{x, y} \in \Z G, $$
where $x,y$ are the colors at the given crossing
($x$ is assigned on the under-arc from which the normal of the over-arc
points, and $y$ is assigned to the over-arc), 
and the product 
$$ ({\cal C}(a_1)^{\epsilon (a_1) } {\cal C}( a_2)^{\epsilon (a_2) }
 \cdots{\cal C}( a_k)^{\epsilon (a_k) } )  \in G_X$$ 
acts on $G$ via 
the quandle module structure. 
The sign $\pm$ in front is determined by whether $r$ is positive $(+)$
or negative $(-)$. The exponent $\epsilon (a_j)$ is $1$ if the arc $\gamma$
crosses the arc $a_j$ against its normal, and is $-1$ otherwise,
for $j=1, \ldots, k$. 
} \end{definition}

\begin{figure}[h]
\begin{center}
\mbox{
\epsfxsize=2in
\epsfbox{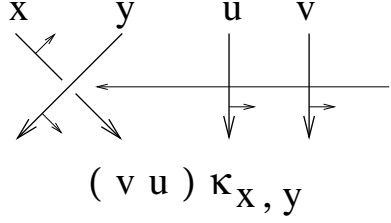} 
}
\end{center}
\caption{Weights on arcs}
\label{weightR} 
\end{figure}

The situation is depicted in Fig.~\ref{weightR}, 
when the arc intersects two arcs with colors $u$ and $v$ in this order.

\begin{figure}[h]
\begin{center}
\mbox{
\epsfxsize=4in
\epsfbox{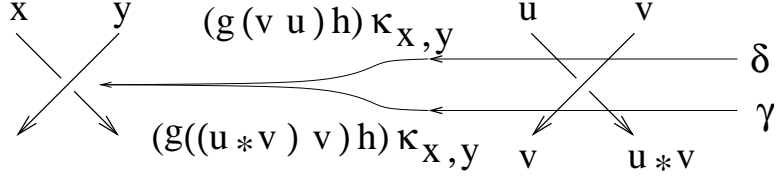} 
}
\end{center}
\caption{When an arc passes a crossing}
\label{pathcross} 
\end{figure}

\begin{lemma} 
The Boltzmann weight,
 does not depend on the choice of the 
arc $\gamma$, so that it will be  denoted by  $B( {\cal C}, r)$. 
\end{lemma} 
{\it Proof.\/} 
It is sufficient to check the changes of the coefficient
$ ({\cal C}(a_1)^{\epsilon (a_1) } {\cal C}( a_2)^{\epsilon (a_2) }
 \cdots{\cal C}( a_k)^{\epsilon (a_k) } )  \in G_X$ 
when the arc is homotoped. 
When a pair of intersection points with the knot diagram 
is canceled or introduced
by a path that  zig-zags, then their colors are inverses of each other,
and are adjacent in the above sequence, so that the product does not 
change. 
It remains  to check what happens when an arc is homotoped through 
each crossing, and the effect of the incident is depicted in 
Fig.~\ref{pathcross}. 

There are two arcs, $\delta$ and $\gamma$, depicted in the figure.
The given crossing $r$ is the left crossing with colors $x$ and $y$. 
The arc $\gamma$ is obtained from $\delta$ by homotoping $\delta$ 
through a crossing with colors $u$ and $v$ as depicted. 
One sees that $B( {\cal C},r, \delta)=(g (vu)h) \kappa_{x,y}$,
where  $g$ and $h$ are sequences of colors  
that $\delta$ intersects before and 
after, respectively, it  intersects  $u$ and $v$.
For the arc $\gamma$, we have  
$B( {\cal C},r, \gamma)=(g ((u*v)v))h) \kappa_{x,y}$, which 
agrees with 
 $B( {\cal C},r, \delta)$. 
Alternative orientations, and signs of crossings follow similarly. 
\qed

\begin{definition} {\rm
The family 
$\Phi_{\kappa} (K) =  
\{ \sum_{r}  B( {\cal C}, r) \}_{{\cal C}}
$
is called the quandle cocycle invariant with respect to 
the (generalized) $2$-cocycle $\kappa$. 
} \end{definition}

The invariant agrees with the quandle cocycle invariant $\Phi_{\phi} (K)$
defined  in \cite{CJKLS}  when the quandle module structure 
of the coefficient group $G$ is trivial, 
and with $\Phi_{\phi} (K)$ defined  in \cite{CES},
modulo the Alexander numbering convention in 
the Boltzmann weight in \cite{CES},
when the coefficient 
group $G$ is 
a $\Z [t, t^{-1}]$-module and the quandle module 
structure is given by $\eta_{x,y}(a)=ta$ and $\tau_{x,y}(b)=(1-t)b$,
for $x, y \in X$, $a, b \in G$.
In the above papers, the state-sum form is used, instead 
of families.
In this section we use families since it is easier to write
for families of vectors.

We note that this definition contains the following data that were
chosen and fixed: $X$, $G$, $\kappa$, and $K$. 

\begin{figure}[h]
\begin{center}
\mbox{
\epsfxsize=4in
\epsfbox{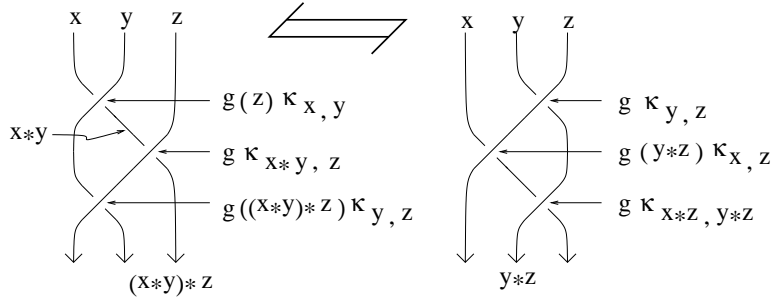} 
}
\end{center}
\caption{Reidemeister type III move and $2$-cocycle condition}
\label{weight3R} 
\end{figure}

\begin{theorem}
The family 
$\Phi_{\kappa} (K) $ does not depend on the choice of a diagram
of a given knot, so that it is a well-defined knot invariant.
\end{theorem}
{\it Proof.\/}
The proof is a routine check of Reidemeister moves, and
it is straightforward for type I and II. 
The type III case is depicted, for one of the orientation choices, 
in Fig.~\ref{weight3R}, where $g$ denotes the sequence of 
colors of arcs that appear before the arc $\gamma$ intersects 
 the crossings in consideration. 
It is seen from the figure that the contribution of the Boltzmann
weights from these three crossings involved do not change before 
and after the move. The other orientation possibilities 
can be checked similarly. 
\qed

\begin{figure}[h]
\begin{center}
\mbox{
\epsfxsize=4in
\epsfbox{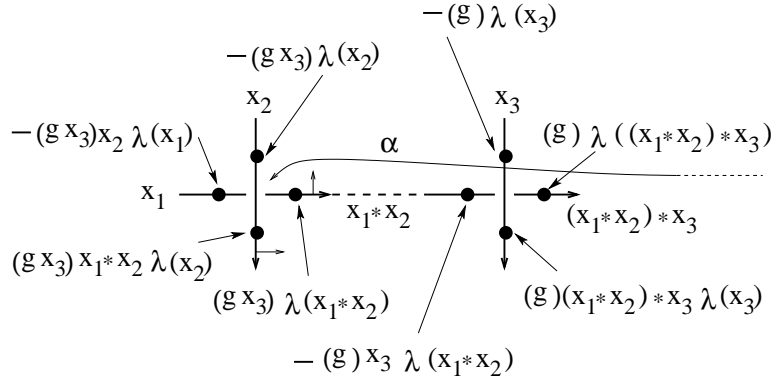} 
}
\end{center}
\caption{Coboundaries cancel}
\label{coboundary} 
\end{figure}

\begin{lemma} \label{coblemma}
If $\kappa=\delta \lambda$ for some $\lambda \in C^1(X;A)$, then
the cocycle invariant $\Phi_\kappa(L)$ is trivial for any 
link $L$.
 Moreover cohomologous cocycles yield the same invariants.
\end{lemma}
{\it Proof.\/} The proof follows the same idea as in \cite{CJKLS}, 
and Proposition~\ref{cobprop}. 
We prove the case when $\delta \lambda =\kappa$,
 the case of cohomologous cocycles follows similarly.
The coboundary of $\lambda$ is given as 
$$\delta \lambda(x,y)=
 -[y \lambda(x) +\lambda(y) -x*y \lambda(y) - \lambda(x*y)].$$ 
Near a crossing with source colors $x,y\in X$, assign $-y\lambda(x) $ 
to the source lower arc, $\lambda(x*y)$ to the target under-arc, 
$-\lambda(y)$ to the over-arc away from which the normal to the
under-arc points, 
and $x*y \lambda(y)$ to the remaining over-arc,
with the action of the sequence of group elements depending on 
the proximity of the arcs to the region at infinity. 
The situation is depicted in Fig.~\ref{coboundary}.
In the figure, the arc from the region at infinity is depicted 
and named $\alpha$, which crosses the arc colored by $x_3$.
The arc $\alpha$ crosses other arcs, and the product of colors 
is denoted by $g$, so that the left crossing has the contribution
 $(gx_3)\kappa_{x_1, x_2}$ to the Boltzmann weight. 
The values $-(gx_3) x_2 \lambda(x_1)$, 
$-(gx_3) \lambda(x_2)$, $(gx_3) x_1 * x_2  \lambda(x_2)$, 
and  $(gx_3) \lambda( x_1 * x_2)$ are assigned to the four arcs
involved at the left crossing. The term  $(gx_3) \lambda( x_1 * x_2)$
cancels with the same term assigned to the right of the horizontal arc
coming from the contribution from the right crossing. 
Thus the assignments to arcs along consecutive crossings cancel. \qed

\begin{proposition} 
In the case when the quandle $2$-cocycle $\kappa$ is expressed
by a group $2$-cocycle as in Proposition~\ref{qmoduleprop}, 
 the invariant $\Phi_{\kappa}(K)$
is trivial for knots $K$ (not necessarily for links).
\end{proposition}
{\it Proof.\/}
Let $K$ denote a knot (not a link). Let $\pi=\pi_1(S^3\setminus K)$ denote 
the  fundamental group  of the complement. 
The fundamental quandle $\pi_Q=\pi_Q(K)$ (see for example \cite{FR})
has $\pi$ as its enveloping group $G_{\pi_Q}=\pi$.
Let $E=A\rtimes_\theta G$ be an extension of a group $G$ by 
an abelian group $A$ twisted by 
a group $2$-cocycle $\theta \in Z^2 _{\rm Group}(G;A)$.  
A quandle coloring of $K$ by ${\mbox{\rm Conj}}(G)$ is induced by a 
quandle homomorphism $f: \pi_Q \rightarrow {\mbox{\rm Conj}}(G)$ 
which naturally lifts to a group homomorphism $f:\pi \rightarrow G$,
denoted   by the same letter. There is an action of $\pi$ on 
$A$ given via the group homomorphism $f$. Explicitly, if $x \in \pi$ and 
$a\in A$, then 
$x\cdot a=s(f(x))i(a)s(f(x))^{-1}$ where $s:G\rightarrow E$ is 
the section that gives rise to the 
description $E=A\rtimes_\theta G$, 
and $i:A\rightarrow E$ is the inclusion. 

The sphere  theorem gives that $S^3\setminus K$ is a $K(\pi, 1)$-space, 
and so $H^2_{\rm Group}(\pi;A)$ is trivial. 
Hence any group $2$-cocycle is a coboundary for $\pi$.
This implies that for generators $\alpha$, $\beta$ of $\pi$ or $\pi_Q$
corresponding to arcs of a given  knot diagram, 
$\theta_{f(\alpha), f(\beta)}=f^{\#}\theta (\alpha, \beta)=
\delta_G \gamma$, and by Lemma~\ref{cobcoblemma}, 
$\kappa_{f(\alpha), f(\beta)}=f^{\#}\delta_Q \gamma$
Then Lemma~\ref{coblemma} implies that the invariant is trivial. 
\qed

\subsection*{Braids and the cocycle invariants}

Let a knot or a link $L$ 
be represented as a closed braid
$\hat{w}$, where $w$ is a $k$-braid word. 
Let $X$ be a quandle and $G$ be a quandle module, keeping the setting 
given at the beginning of this section.

\begin{definition} {\rm 
Let $\vec{x}=(x_1, \ldots, x_k)\in X^k$ for some positive integer $k$.
Define a {\it weighted sum} on $G^k$ with respect to $\vec{x}$
 by
$$ \mbox{WS}_{\vec{x}}(\vec{a})= 
\sum_{i=1}^k u_i a_i = 
(x_k \cdots x_2  ) a_1 +   (x_k \cdots x_3) a_2 
+ \cdots +  x_{k} a_{k-1} +  a_k,$$
where $\vec{a} =(a_1, \ldots, a_k) \in G^k$ 
and $u_i=x_k \cdots x_{i+1}$ for $i=1, \ldots, k-1$, and $u_k=1$.
} \end{definition}

Let $w$ be a $k$-braid word,
$\vec{x}=(x_1, \ldots , x_k) \in X^k$ be a vector of colors assigned
to the bottom strings of $w$, and ${\cal C}$ be the unique coloring of $w$ 
by $X$ determined by $\vec{x}$. 
Let the vector $\vec{y}=(y_1, \ldots y_k) \in X^k$ denote the  color vector
on the top strings  that is induced by the vector $\vec{x}=(x_1, \ldots x_k) \in X^k$ which is assigned at the bottom. 
Fix this coloring ${\cal C}$.
Let $\alpha^0_{x,y}= \eta_{x,y} + \tau_{x,y} $
be a quandle $2$-cocycle with coefficient group $G$
with $\kappa=0$, 
where $G$ is an $X$-module.
Then any  color vector $\vec{a}=(a_1, \ldots, a_k) \in G^k$ 
at the bottom strings uniquely determines a color vector 
$\vec{b}=(b_1, \ldots, b_k) \in G^k$ at the top strings of $w$
with respect to $\alpha$, that is, $\vec{b}=M(w, \vec{x}) \cdot \vec{a}$.
Recall that in this section the quandle module structure is given by 
the $G_X$-module structure.

\begin{lemma} \label{wslemma}
Let  $\alpha^0_{x,y}= \eta_{x,y} + \tau_{x,y} $ be a quandle $2$-cocycle 
with $\kappa=0$, and $\vec{b}=M(w, \vec{x}) \cdot \vec{a}$, as above.
Then  we have 
$\mbox{WS}_{\vec{x}}(\vec{a})=\mbox{WS}_{\vec{y}}(\vec{b})$.
\end{lemma}
{\it Proof.\/}
It is sufficient to prove the statement for a generator
and its inverse. The inverse case is similar, so we compute 
the case when $w=\sigma_i$ for some $i$, $1 \leq i < k$. 
The only difference in the weighted sum, in this case, 
is the $i\/$th and $(i+1)$st terms. 
Let $u=x_k \cdots x_{i+2}$. 
For the bottom colors, the $i\/$th and $(i+1)$st terms are 
$$ \mbox{WS}_{\vec{x}}(\vec{a})=
\cdots + u x_{i+1}  a_i + u   a_{i+1} + \cdots .$$
On the other hand, one computes
$$
 \mbox{WS}_{\vec{y}}(\vec{b}) = 
\cdots + u (x_i * x_{i+1}) a_{i+1} 
+ u 
[ x_{i+1} a_i + (1-  x_i * x_{i+1} ) a_{i+1} ] + \cdots $$
which agrees with 
the above.
\qed

\begin{theorem} \label{wsthm}
Let $w$ be a $k$-braid word with crossings $\rho_{\ell}$,
$\ell=1, \ldots, h$, and $\vec{x}$ a bottom color vector
by a quandle $X$. 
Let  $\alpha_{x,y}= \eta_{x,y} + \tau_{x,y} + \kappa_{x,y} $ 
be a quandle $2$-cocycle with coefficient 
group $G$, a $G_X$-module, 
and $\vec{b}=M(w, \vec{x}) \cdot \vec{a}$.
Here, the map $M$ corresponds to $\alpha$ with possibly non-zero $\kappa$.
Then we have 
$$ \mbox{WS}_{\vec{y}}(\vec{b}) - \mbox{WS}_{\vec{x}}(\vec{a})
= \sum_{\ell=1}^{h}  B({\cal C},\rho_{\ell}). $$
\end{theorem} 
{\it Proof.\/} 
Let $M^0(w, \vec{x})$ denote the map corresponding to 
the $2$-cocycle  $\alpha^0_{x,y}= \eta_{x,y} + \tau_{x,y}$,
which is obtained from $\alpha$ by setting $\kappa=0$.
The theorem follows from induction once we prove it for 
a braid generator and its inverse, and we show the case
$w=\sigma_i$, as the inverse case is similar.
Then we compute
\begin{eqnarray*} 
\mbox{WS}_{\vec{y}}(\vec{b})=\mbox{WS}_{\vec{y}}( M(\sigma_i, \vec{x} )\cdot  \vec{a} ) &=& 
\mbox{WS}_{\vec{y}}( M^0(\sigma_i, \vec{x} )\cdot \vec{a} )
+ B({\cal C}, \sigma_i) \\
&=& \mbox{WS}_{\vec{x}}( \vec{a} ) - B({\cal C}, \sigma_i),
\end{eqnarray*}
where the first equality follows from the definitions,
and the second equality follows from Lemma~\ref{wslemma}.
Note that the braid generator represents a negative crossing
in the definition of the quandle module invariant, and a positive
crossing in the cocycle invariant, so that there is a negative sign 
for the weight $B({\cal C}, \sigma_i)$. 
\qed

\begin{theorem} \label{colorextthm}
If a  coloring ${\cal C}$ of $L$ by $X$  extends
 to a coloring of $L$ by the extension $E= G \times_{\alpha} X$,
then the coloring contributes 
a trivial term (i.e., an integer in $\Z G$)
to  the generalized cocycle invariant $\Phi_{\kappa} (L) $.
\end{theorem}
{\it Proof.\/}
A given coloring agrees  on the bottom and top strings, so that 
$\vec{x}=\vec{y}$. 
If the given coloring extends to $E$, we have 
$ \vec{b}= M(\sigma_i, \vec{x} )\cdot  \vec{a}=  \vec{a}$,
and in particular, 
$\mbox{WS}_{\vec{y}}( \vec{b})= \mbox{WS}_{\vec{x}}( \vec{a} )$.
Then this theorem follows from Theorem~\ref{wsthm}.
\qed

5Thus the invariant  $\tilde{\Phi}(X, \alpha\ ; L)$ measures exactly

\begin{remark} {\rm 
In Livingston \cite{Liv1995}, the following situation 
is examined. Suppose that $E$ is a split extension of a group $G$ by 
an abelian group $A$, 
so that $0\rightarrow A \rightarrow E \rightarrow G \rightarrow 1$ is
 a split exact 
sequence. 
Let $\rho :\pi=\pi_1(S^3\setminus K) \rightarrow G$ denote a homomorphism. 
Then since there is an action of $G$ on $A$ via conjugation in $E$, 
there is a corresponding action of $\pi$ on $A$ given via $\rho$. 
Livingston examines when there is a lift of $\rho$ to 
$\tilde{\rho}: \pi \rightarrow E$ thereby generalizing Perko's theorem 
that any homomorphism $\rho: \pi \rightarrow \Sigma_3$ lifts to 
a homomorphism $\tilde{\rho}: \pi \rightarrow \Sigma_4$. 
In this situation, the permutation group $\Sigma_3$ acts on 
the Klein $4$ group 
$\Z_2 \times \Z_2 = \{ (1), (12)(34), (13)(24), (14)(23) \}$ 
via conjugation and the obvious section $s:\Sigma_3 \rightarrow \Sigma_4$. 
Livingston shows that  there is a one-to-one correspondence between 
$A$-conjugacy classes of lifts of $\rho$  and elements in $H^1(\pi, \{A\})$ 
where the  coefficients
$\{A\}$ are twisted by 
the action of $\pi$  on $A$.

Consider a quandle coloring of a knot $K$ by a quandle 
${\mbox{\rm Conj}}(G)$. 
Such a quandle coloring is a homomorphism from the fundamental quandle, 
$\pi_Q(K)$, of $K$,
to ${\mbox{\rm Conj}}(G)$. The fundamental group $\pi$  can be thought of as 
$\pi=G_{\pi_Q(K)}$. Thus a quandle coloring is a homomorphism $\rho$ as above. 

Thus, we have the following  cohomological characterization of lifting 
colorings. If a coloring lifts, then the coloring 
contributes a constant 
term to the cocycle invariant
by Theorem~\ref{colorextthm} above, 
 and it corresponds to a $1$-dimensional cohomology class 
of $H^1(\pi, \{A\})$. 
} \end{remark}

\subsection*{Computations}

In this section we present a computational method using closed braid form.
{\it We take a positive crossing as a positive generator, in this section.} 
Let $w=w_1 \cdots w_h$ be a braid word, 
where $w_s=\sigma_{j(s)}^{\epsilon(s)}$ is a standard generator or its inverse
for each $s=1, \cdots, h$. The braids are oriented downward and
the normal points 
to the right. Then the target region is to the right 
of each crossing. 

\begin{figure}[h]
\begin{center}
\mbox{
\epsfxsize=3in
\epsfbox{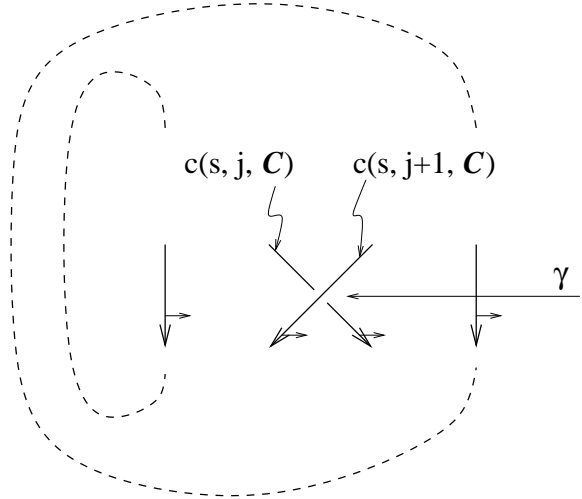} 
}
\end{center}
\caption{Colors and arcs for braids }
\label{braidconvention} 
\end{figure}

Let $X$ be a quandle and ${\mathcal C}$ be a coloring of $\hat{w}$ 
by $X$. Let $c(s, i,  {\mathcal C})$ be the color of 
the $i\/$th string from the left, 
immediately 
above 
the $s\/$th 
crossing ($w_s$) for ${\mathcal C}$. 
Take the right-most region as the region at infinity,
and take the closure of a given braid to the left, 
as depicted in Fig.~\ref{braidconvention}. 
As an arc  $\gamma$ to the target region, take an arc that goes 
horizontally 
from the right to left to the crossing, see Fig.~\ref{braidconvention}.
Then the Boltzmann weight at the  $s$-th crossing is 
given by 
$$ B(  {\mathcal C}, w_s)
= 
(c(s, n, {\mathcal C})  \cdots   c(s, j+2, {\mathcal C}))
\kappa_{c(s, j,  {\mathcal C}),  c(s, j+1,  {\mathcal C})}
$$ 
for a positive crossing (when $\epsilon(s)=1$), 
and the $2$-cocycle evaluation is replaced by 
$\kappa_{c(s+1, j,  {\mathcal C}),  c(s, j,  {\mathcal C})}$
for a negative crossing, and with a negative sign in front.
Note that the expression in front of $\kappa$ represents 
the action of this group element on the coefficient group, so that 
in the case of wreath product, for example, this can be written 
in terms of matrices.
Also, if $j=n-1$, then the  expression in front is understood to be
empty. 
This formula can be used to evaluate the invariant for knots 
in closed braid form, and can be implemented 
in a computer.

\begin{example} \label{r3classicalex} {\rm   
We implemented the above formula for $R_3$ in {\it Maple} and 
obtained the following results,
confirmed also by {\it Mathematica}. 
The action is 
in the wreath product; thus $R_3$ acts  
on $(\Z_q)^3$ by 
transpositions of factors.
We also represent elements of $R_3$ by $\{ 1,2,3 \}$ 
with the following correspondence in this example:
 $1=(2\ 3)$,  $2=(1\ 3)$ and $3=(1\ 2)$.

For $q=0$, any $2$-cocycle gave trivial invariant 
(each coloring contributes the zero vector, and the family
of vectors is a family of zero vectors),
for all $3$-colorable knots in the knot table up to $9$ crossings.
Thus we conjecture that $q=0$ gives rise to the trivial invariant.

Let $h(i,j)={}^T( f_1(i,j),  f_2(i,j),  f_3(i,j)) $
denote a vector valued $2$-cochain. Then  
for $q=3$, the following defines a $2$-cocycle:
\begin{eqnarray*}
 f_1(3, 1)\ =\ f_3(3,1)\ =\  f_1(2,3)\ =\  f_1(2,1)\ =\  f_3(2,1) &= & 1, \\
f_2(3, 1)\  =\ f_2(3,2)\ =\ f_2(1,3)\ =\  f_1(1,2)\ =\   f_2(1,2) &= & 2
  \end{eqnarray*} 
and all the other values are zeros. 
Then we obtain the following results.
\begin{itemize}
\setlength{\itemsep}{-3pt}
\item
$\Phi_{\kappa}(K)=\{ \sqcup_9 (0,0,0) \}$ for 
$K=  6_1,  8_{10}, 8_{11}, 8_{20}, 9_1, 9_6, 9_{23}, 9_{24}$,
where $\{ \sqcup_9 (0,0,0) \}$ represents the 
 family consisting of $9$ copies of  $(0,0,0)$
(similar notations are used below).
\item
$\Phi_{\kappa}(K)=\{ \sqcup_3 (0,0,0),\sqcup_6 (1,1,1)  \}$ for 
$K=3_1, 7_4, 7_7, 9_{10}, 9_{38}$.
\item
$\Phi_{\kappa}(K)=\{ \sqcup_3(0,0,0), \sqcup_6 (2,2,2) \} $ for 
$K=8_5, 8_{15}, 8_{19}, 8_{21}, 9_2, 9_4, 9_{11}, 9_{15}, 9_{16}, 9_{17}, 9_{28}9_{29}, 9_{34}, 9_{40}$.
\item
$\Phi_{\kappa}(K)=\{  \sqcup_9 (0,0,0), \sqcup_{18} (1,1,1) \} $ for
 $K= 9_{35}, 9_{47}, 9_{48}$.
\item
$\Phi_{\kappa}(K)=\{  \sqcup_3 (0,0,0), \sqcup_{12} (1,1,1), \sqcup_{12} (2,2,2) \} $ for $K=8_{18}$.
\item
$\Phi_{\kappa}(K)=\{  \sqcup_{15} (0,0,0), \sqcup_6(1,1,1), \sqcup_6 (2,2,2)\} $ for $K=9_{37}, 9_{46}$.
\end{itemize} 
In the original cohomology theory with trivial action on the coefficient,
the second cohomology group $H^2_{\rm Q}(R_3; G)$ was trivial 
for any coefficient group $G$ \cite{CJKLS}, so that $R_3$ gave 
rise to a trivial cocycle invariant. This example shows that, 
with the wreath product action, $R_3$ has non-trivial second cohomology group, 
and gives rise to a non-trivial cocycle invariant, showing that the theories 
with actions on coefficients are strictly more general and stronger than 
the original case.
} \end{example}

\begin{remark} {\rm
With the same cocycle as the preceding example, we computed
the invariants for the mirror images with the same orientations.
The results are such that the values $1$ and $2$ are exchanged
in all values
(thus we conjecture that this is the case in general,
at least for $R_3$). 
For example, the mirror image of  the trefoil has 
$\Phi_{\kappa}(K)=\{ \sqcup_3 (0,0,0),\sqcup_6 (2,2,2)  \}$
as its invariant. 
Hence those with asymmetric values of the invariant 
are proven to be non-amphicheiral by this invariant. 
Specifically, $22$ knots among $33$ are detected to be non-amphicheiral.
 
} \end{remark}

\section{Invariants for knotted surfaces}\label{knottedsfcesec} 

\subsection*{Definitions}

The cocycle invariants are defined for knotted surfaces in $4$-space 
 in  exactly the same manner as in Section~\ref{cocyinvsec} as follows.
Let $X$ be a finite quandle and $\Z (X)$ be its quandle algebra
with generators $\{ \eta _{x,y}^{\pm 1}  \}_{ x, y \in X} $
and $\{ \tau _{x,y}  \}_{ x, y \in X} $.
Let $G$ be 
an abelian group that is a $G_X$-module. Recall that 
this induces a 
$\Z (X)$-module structure  
given by $\eta_{x,y} g = y g$ and $\tau_{x,y}(g) = (1 - x*y)g$
 for $g \in G$ and 
$x,y \in X$. 
Let $\kappa_{x,y,z}$ be a generalized quandle $3$-cocycle of $X$ with 
the coefficient group $G$. 
Thus the generalized $3$-cocycle condition, in this setting, is 
written as 
\begin{eqnarray*}
\lefteqn{
 w \kappa_{x,y,z} + \kappa_{x*z, y*z, w} + ((y*z)*w) \kappa_{x,z,w}
+ \kappa_{y,z,w} } \\
&= & ((x*(y*z))*w)  \kappa_{y,z,w}  + \kappa_{x*y, z, w} 
+ (z*w) \kappa_{x,y,w} + \kappa_{x*w, y*w, z*w}.
\end{eqnarray*}
We require further that 
$\kappa_{x,x,y}=\kappa_{x,y,y}=0$.
A cocycle invariant of knotted surfaces
 will be defined using such a $3$-cocycle. 

A knotted surface diagram $K$ is given in  $3$-space. 
We assume the surface is oriented and 
use orientation normals to indicate the orientation. 
In a neighborhood of each triple point, there are eight regions 
that are separated by the sheets of the surface since
 the triple point looks like the intersection of the 
$3$-coordinate planes in some parametrization.
The region into which 
all normals point
is called the {\it target} region.  
Let $\gamma$ be an arc from the region at infinity of the $3$-space
 to the target region of a given triple point $r$. Assume  that 
$\gamma$
intersects
$K$ transversely in a finitely many points thereby missing 
double point curves, branch points, and triple points. 
Let $a_i$, $i=1, \ldots, k$, in this order, be the sheets of $K$
that intersect $\gamma$ from the region at infinity to the triple point $r$.
Let ${\cal C}$ be a coloring of $K$ by a fixed finite quandle $X$.

\begin{definition}
{\rm The {\em Boltzmann weight} $B( {\cal C}, r,\gamma)$
for the triple point $r$, for a coloring ${\cal C}$, with respect to $\gamma$,
is defined by 
$$ B( {\cal C}, r,\gamma)=
 \pm  ({\cal C}(a_1)^{\epsilon (a_1) } {\cal C}( a_2)^{\epsilon (a_2) }
 \cdots{\cal C}( a_k)^{\epsilon (a_k) } ) 
 \kappa_{x, y, z} \in \Z G, $$
where $x,y,z$ are the colors at the given triple point $r$
($x$ is assigned on the bottom sheet from which the normals of the
middle and top sheets
point, and $y$ is assigned to the middle sheet from which the normal
of the top sheet points, and $z$ is assigned to the top sheet). 
The sign $\pm$ in front is determined by whether $r$ is positive $(+)$
or negative $(-)$. The exponent $\epsilon (a_j)$ is $1$ is the arc $\gamma$
crosses the arc $a_j$ against its normal, and is $-1$ otherwise,
for $j=1, \ldots, k$. 
} \end{definition}

\begin{lemma} 
The Boltzmann weight does not depend on the choice of the 
arc $\gamma$, so that it will be  denoted by  $B( {\cal C}, r )$. 
\end{lemma} 
{\it Proof.\/} 
It is sufficient to check what happens when an arc is homotoped through 
double point curves,  and  
Fig.~\ref{pathcross} can be 
regarded as a cross-sectional view of such a homotopy,
and the same computation as in the classical case holds. 
\qed

\begin{definition} {\rm
The family 
$\Phi_{\kappa} (K) = 
\{ \sum_{r} B( {\cal C}, r) \}_{{\cal C} \in {\rm Col}_{\rm X} (K) } $
is called the quandle cocycle invariant with respect to 
the (generalized) $3$-cocycle $\kappa$. 
} \end{definition}

\begin{figure}[h]
\begin{center}
\mbox{
\epsfxsize=5.5in
\epsfbox{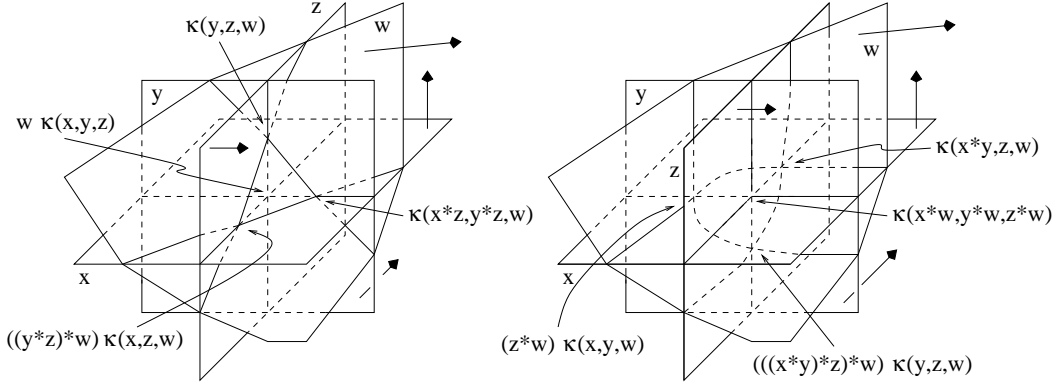} 
}
\end{center}
\caption{Contributions for the tetrahedral move }
\label{tetranu} 
\end{figure}

\begin{theorem}
The family 
$\Phi_{\kappa} (K) $ does not depend on the choice of a diagram
of a given knotted surface, so that it is a well-defined knot invariant.
\end{theorem}
{\it Proof.\/}
The proof is a routine check of Roseman moves, analogues of
Reidemeister moves. 
In particular, the analogue of the type III Reidemeister move 
is called the tetrahedral move, and is a generic plane passing through
the origin that is the triple point formed by coordinate planes.
In Fig.~\ref{tetranu}, such a move is depicted.
A choice of normal vectors and quandle colorings are also depicted
in the figure.
The sheet labeled $w$ is the top sheet, and the next highest sheet
is labeled by $z$, the bottom is $x$. 
The region at infinity is chosen (for simplicity) 
to be the region at the top right, into which
all normals point. The Boltzmann weight at each triple point is also 
indicated. The sum of the weights for the LHS and RHS are exactly
those for the $3$-cocycle condition. 
Other choices for orientations, and other moves, are checked similarly.
In particular, by assuming that the cocycle $\kappa$ satisfies the 
quandle condition $\kappa_{x,x,y}=\kappa_{x,y,y}=0$, 
we ensure that the quantity is invariant under the
Roseman move in which a branch point passes through another sheet.
\qed

\begin{figure}[h]
\begin{center}
\mbox{
\epsfxsize=3in
\epsfbox{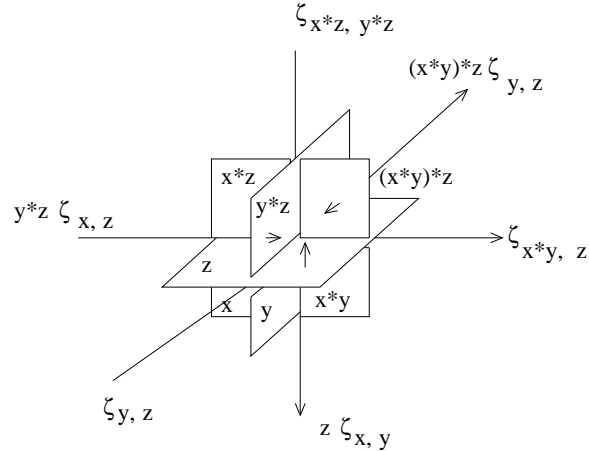} 
}
\end{center}
\caption{Coboundary terms distributed  }
\label{cocycleattrip} 
\end{figure}

A proof similar to that of Lemma~\ref{coblemma} 
implies the following, where the distribution of coboundary terms 
is depicted in Fig.~\ref{cocycleattrip}.

\begin{lemma}
If $\kappa=\delta \zeta$ for some $\zeta \in C^2(X;A)$, then
the cocycle invariant $\Phi_\kappa(F)$ is trivial for any 
knotted surface $F$.
 Moreover cohomologous cocycles yield the same invariants.
\end{lemma}

\subsection*{Computations}

\begin{figure}[h]
\begin{center}
\mbox{
\epsfxsize=3.5in
\epsfbox{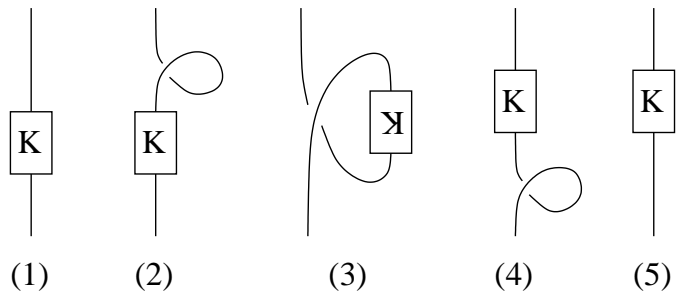} 
}
\end{center}
\caption{A movie of twist spinning}
\label{twistmovienu} 
\end{figure}

We develop a computational method, 
based on Satoh's method \cite{Satoh}, of computing this cocycle invariant 
for twist-spun knots using 
the  closed braid form.
First we review Satoh's method.
A movie description of one full twist of a classical knot $K$
is depicted in Fig.~\ref{twistmovienu} from (1) through (5). 
Strictly speaking, $K$ is a tangle with two end points, and 
the tangle is twisted 
about an axis containing these end points.
For the twist-spun trefoil, a diagram of the trefoil (with two end points)
goes in the place
of $K$.

\begin{figure}[htb]
\begin{center}
\mbox{
\epsfxsize=3.5in
\epsfbox{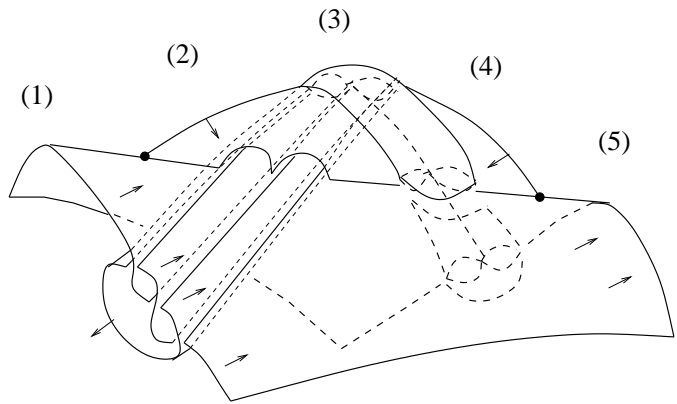} 
}
\end{center}
\caption{A part of a diagram of twist-spun trefoil}
\label{spuntrediagnu} 
\end{figure}

It is seen from this movie that branch points 
appear between (1) and (2), and (4) and (5), when type I Reidemeister
moves occur in the movie. 
Triple points appear when $K$ goes over an arc between 
(2) and (3), and goes under it between (3) and (4). 
Each triple point corresponds to a crossing of the diagram $K$.  
 This movie will provide a broken surface 
diagram by taking the continuous traces of the movie,
which is depicted in  Fig.~\ref{spuntrediagnu}.
Horizontal cross sections of Fig.~\ref{spuntrediagnu}
 correspond to (1) through (5)
as indicated at the top of the figure. Thus  between 
(1) and (2) there is a branch point, for example.
A choice of normal vectors is also depicted by short arrows.

\begin{figure}[htb]
\begin{center}
\mbox{
\epsfxsize=1.2in
\epsfbox{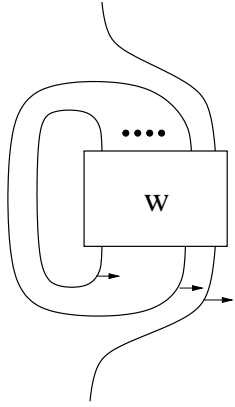} 
}
\end{center}
\caption{A tangle construction from a braid}
\label{braidtangle} 
\end{figure}

Now we  apply this method to closed braid form.
Let $w$ be a diagram corresponding to a braid word
of $n$ strings. Then construct a tangle with two end points 
at top and bottom, by stretching the left-most end points 
and closing all the other end points of the braid, as 
depicted in Fig.~\ref{braidtangle}.
Let $K$ be this tangle, as well as the corresponding knot
and knot diagram.
A choice of normal vectors are also depicted.
Let ${\rm Tw}^{\ell}(K)$ be the diagram of the $\ell$-twist 
spun of $K$ obtained
by Satoh's method.

\begin{figure}[htb]
\begin{center}
\mbox{
\epsfxsize=3.5in
\epsfbox{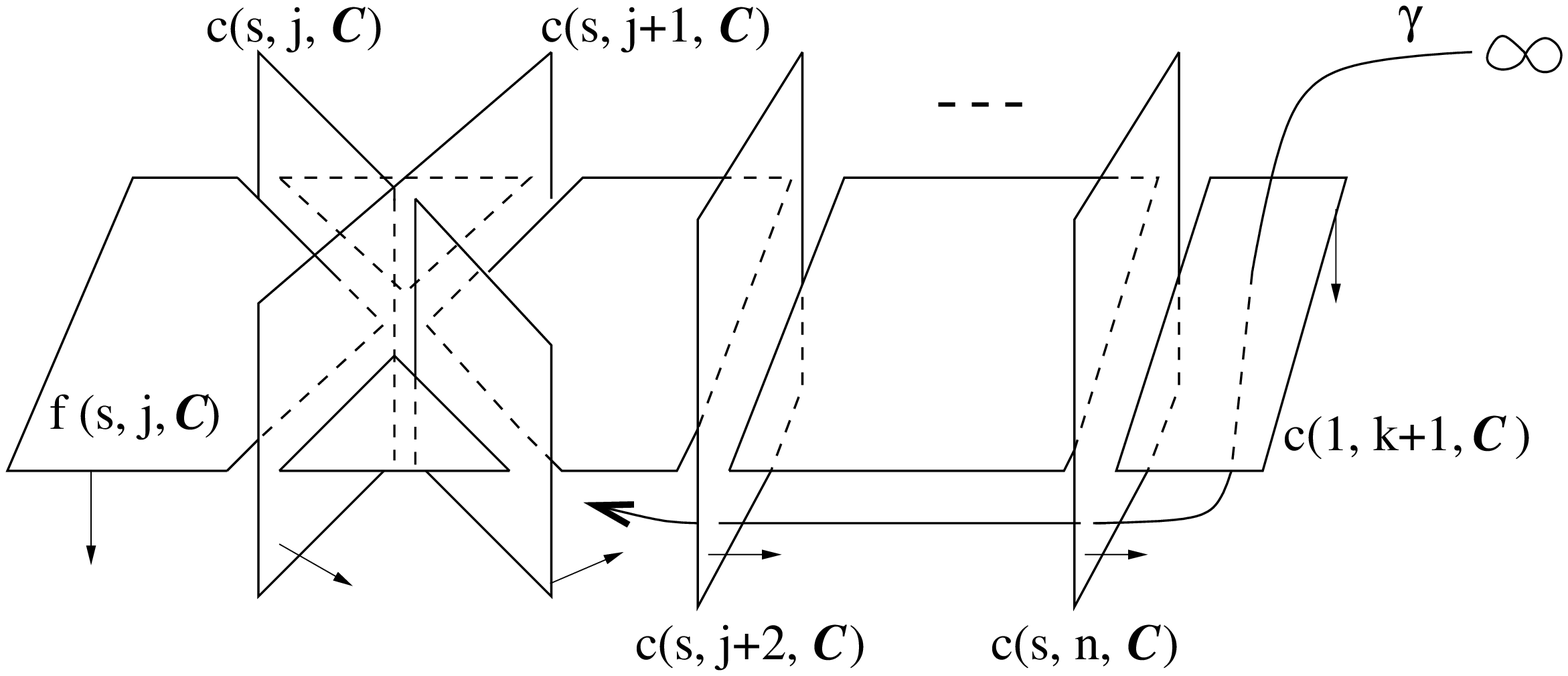} 
}
\end{center}
\caption{The weight at a triple point}
\label{triplewt} 
\end{figure}

Let $w$ (as a braid word, written as the same letter as the diagram), 
be $w_1, \ldots, w_h$, where each $w_s$ is a standard generator
or its inverse, $\sigma_{j(s)}^{\epsilon(s)}$. 
In this section we use the positive crossing 
as a standard generator. Recall that 
each crossing gives rise to a triple point 
in the diagram of  a twist-spun knot, when Satoh's method 
is applied. In Fig.~\ref{triplewt}, the triple point
corresponding to $w_s=\sigma_j$ is depicted.
This figure represents a triple point $T_1^{l}(s)$ that is formed when 
the crossing $w_s$ goes through a sheet between steps (2) and (3) 
in Fig.~\ref{spuntrediagnu} for  $j=j(s)$ and $\epsilon(w_s)>0$
(i.e., $w_s$ is the $j$-th braid generator $\sigma_j$, 
and the crossing is  positive). 
There is another triple point $T_1^{r}(s)$ formed by the same crossing between
steps (3) and (4). The superscripts $l$ and $r$ represents that 
they appear in the left and right of  Fig.~\ref{spuntrediagnu},
respectively. It is seen that $T_1^{l}(s)$
and $T_1^{r}(s)$ are negative and positive triple points with the right-hand
rule, respectively, with respect to the normal 
vectors specified in Figs.~\ref{spuntrediagnu}
and \ref{braidtangle}. 
Then there are a pair of  triple points $T_u^{l}(s)$ and $T_u^{r}(s)$
for the left and right, respectively, for the $u$-th twist.

Let a coloring ${\mathcal C}$ by a quandle $X$ of 
the diagram ${\rm Tw}^{\ell}(K)$ be given.
At the triple point $T_1^{l}(s)$, the triple of colors that 
contribute to the cocycle invariant  are depicted 
in Fig.~\ref{triplewt} 
as 
$f(s, j, {\mathcal C})$, $c(s, j, {\mathcal C})$ and $c(s, j+1, {\mathcal C})$
for the top, middle, and bottom sheet, respectively. 
Here $f$ represents 
a color of the face left to the crossing $w_s$,
and $j$ in $c(s, j, {\mathcal C})$ represents 
the  $j$-th string from 
the left in the braid $w$ (as $w_s=\sigma_j$). 
The color of the horizontal sheet in Fig.~\ref{triplewt} at 
the right-most region outside of the braid is the same as the color 
$c(1,k, {\mathcal C})$ 
of the $k$-th string of the braid $w$ at the top and bottom,
since the tangle goes through the $k$-th string, see Figs.~\ref{braidtangle}
and \ref{twistmovienu}.
Thus we obtain
$$ f(s, j, {\mathcal C}) = ( \cdots ( c(1, k, {\mathcal C}) \bar{*}
c(s, k, {\mathcal C}) ) \bar{*} c(s, k-1, {\mathcal C}) )
 \bar{*}  \cdots )  \bar{*}
c(s, j, {\mathcal C}) ) \cdots ) ,
$$
where for any $b, c\in X$ the unique element $a \in X$ 
with $a*b=c$ is denoted by $a=c \bar{*} b$.  
{}From Figs.~\ref{spuntrediagnu} and \ref{triplewt}, 
it is seen that the target region 
of $T_1^l(s)$ is to the bottom right.
Thus from the region at infinity, we choose an arc $\gamma$ 
as depicted in Fig.~\ref{triplewt}, that goes through
 the  sheet with color $c(1, k, {\cal C})$, then hits
 all the sheets corresponding to the $n$-th through $(j+2)$-th 
braid strings, to the target region.
Only the first sheet is oriented coherently with the direction of
the arc we chose. Hence the sequence of quandle elements 
in the Boltzmann  weight that 
corresponds to the arc $\gamma$ is 
$$ c(1, k, {\mathcal C})^{-1} 
 c(s, n, {\mathcal C}) c(s, n-1, {\mathcal C})  \cdots 
 c(s, j+2, {\mathcal C}). $$
Similar arcs can be chosen for $T_1^r(s)$ as well,
giving the same sequence. 
Thus we obtain the Boltzmann weights as follows. 
$$
  B(  {\mathcal C}, T_1^l(s))
= \left\{ 
\begin{array}{ll}
- ( c(1, k, {\mathcal C})^{-1} 
 c(s, n, {\mathcal C})  \cdots  c(s, j+2, {\mathcal C}) ) 
\kappa_{f(s, j, {\mathcal C}),\
c(s, j,  {\mathcal C}),\  c(s, j+1,  {\mathcal C})} & \quad
\mbox{\rm if}\quad  \epsilon (T_1^l(s)) < 0 
\\
 c(1, k, {\mathcal C})^{-1}
 c(s, n, {\mathcal C})  \cdots  c(s, j+2, {\mathcal C}) 
\kappa_{f(s, j, {\mathcal C}),\ c(s+1, j, {\mathcal C}),\ c(s, j, {\mathcal C})
} 
 & \quad
\mbox{\rm if}\quad \epsilon (T_1^l(s)) > 0 , 
\end{array}
\right.
$$ 
$$ 
  B(  {\mathcal C}, T_1^r(s))
= \left\{ 
\begin{array}{ll}
 c(1, k, {\mathcal C})^{-1} 
 c(s, n, {\mathcal C})  \cdots  c(s, j+2, {\mathcal C}) 
\kappa_{c(s, j, {\mathcal C}),\ c(s, j+1, {\mathcal C}),\
c(1, k, {\mathcal C})} & \quad
\mbox{\rm if} \quad \epsilon (T_1^r(s)) > 0 
\\
-( c(1, k, {\mathcal C})^{-1}
 c(s, n, {\mathcal C})  \cdots  c(s, j+2, {\mathcal C}) )
\kappa_{c(s+1, j, {\mathcal C}),\ c(s, j, {\mathcal C}),\
c(1, k, {\mathcal C})} 
 & \quad
\mbox{\rm if} \quad \epsilon  (T_1^r(s)) < 0 .
\end{array}
\right.
$$

For the $u$-th twist, the cocycle evaluation is replaced by 
$$\kappa_{f(s, j, {\mathcal C}) * c(1,k,{\mathcal C})^u,\ 
c(s, j,  {\mathcal C}) * c(1,k,{\mathcal C})^u,\  
 c(s, j+1,  {\mathcal C}) * c(1,k,{\mathcal C})^u}$$ 
for a positive triple point $T^l_u(s)$, and similar changes 
(multiplication by $ * c(1,k,{\mathcal C})^u$ for every term)
are made for 
the other triple points accordingly, where 
$x*y^{\ell}$ denotes the $\ell$-fold product
$( \cdots (x*y)*y) \cdots )*y$.

With these weights the cocycle invariant is computed by
$$\Phi_{\kappa}({\rm Tw}^{\ell}(K) )
=\left\{ \sum_{u=1}^{\ell} \sum_{s=1}^h \left[
 B(  {\mathcal C}, T_u^l(s)) + B(  {\mathcal C}, T_u^r(s))\right] 
\right\}_{\mathcal C}
 .$$
This formula, again, can be implemented in computer calculations.

\begin{table}[h] 
\begin{center}
{\begin{tabular}{||l||l||}\hline \hline
Knot $K$ & $\Phi_{\kappa}({\rm Tw}^2 (K) )$
\\ \hline \hline 
$3_1$ & $\sqcup_9 (0,0,0). $ 
\\ \hline 
$6_1$ & $\sqcup_3 (0,0,0),\; 
(-q_2, 0, q_2),\;  (q_2, q_1,  - q_1-q_2),\; 
     (q_1, 0, -q_1),\; $
\\ \hline 
 & $  (0, -q_1, q_1),\; 
      ( - q_1+q_2,  q_1 -q_2, 0),\;  (-q_2,  - q_1+q_2, q_1).$  
\\ \hline 
$7_4$  & $ \sqcup_3 (0,0,0),\; 
(-q_2, -2 q_1, 2 q_1 + q_2),\; 
  (-q_1 + q_2, q_1, -q_2),\; $
\\ \hline 
 & $ (-q_1, 2 q_1, -q_1)\; ,
 (-q_1, -q_1, 2 q_1),\;  (q_2, -q_2, 0),\; 
     (-q_2, q_2, 0) . $ 
\\ \hline 
$7_7$    &  $ \sqcup_3 (0,0,0),\;
\sqcup_2 (q_1, 0, -q_1),\; \sqcup_2 (0, -q_1, q_1),\; 
\sqcup_2  (-q_1, q_1, 0).  $ 
\\ \hline 
$8_5$  &  $ \sqcup_9 (0,0,0). $ 
\\ \hline 
$8_{10}$ &  $ \sqcup_9 (0,0,0). $ 
\\ \hline
$8_{11}$   &  $ \sqcup_3 (0,0,0),\; 
 (q_2, 0, -q_2),\;  (-q_2, -q_1,   q_1+q_2),\; 
     (-q_1, 0, q_1),\; $
\\ \hline 
 & $ (0, q_1, -q_1),\; 
  (q_1 - q_2,  - q_1+q_2, 0),\;  (q_2, q_1 - q_2, -q_1). $  
\\ \hline 
$8_{15}$   &  $ \sqcup_3 (0,0,0),\;  
 \sqcup_3 (0, -q_1, q_1),\; \sqcup_3 (-q_1, q_1, 0). $  
\\ \hline 
$8_{18}$    &  $ \sqcup_9 (0,0,0),\;
 \sqcup_6 (q_1, 0, -q_1),\;
 \sqcup_6  (-q_1, q_1, 0),\;
 \sqcup_6 (0, -q_1, q_1).$
\\ \hline
$8_{19}$    &   $ \sqcup_9 (0,0,0). $ 
\\ \hline 
$8_{20} $   &  $ \sqcup_9 (0,0,0).  $  
\\ \hline 
$8_{21}$    &   $ \sqcup_9 (0,0,0). $  
\\ \hline 
\end{tabular} } \end{center}
\caption{A table of cocycle invariants for twist spun knots, part I}
\label{twistspintable1}
\end{table}

\begin{table}[h] 
\begin{center}
{\begin{tabular}{||l||l||}\hline \hline
Knot $K$ & $\Phi_{\kappa}({\rm Tw}^2 (K) )$
\\ \hline \hline 
$9_1$   &    $ \sqcup_9 (0,0,0). $ 
\\ \hline
$9_{2}$    &   $ \sqcup_3 (0,0,0),\; 
\sqcup_2 (0, -q_1, q_1),\;  \sqcup_2 (-q_1, q_1, 0),\;
 \sqcup_2 (q_1, 0, -q_1).  $  
\\ \hline
$9_4$     &   $ \sqcup_3 (0,0,0),\; 
\sqcup_2  (0, q_1, -q_1),\;
\sqcup_2  (q_1, -q_1, 0),\;  \sqcup_2  (q_1, -q_1, 0). $  
\\ \hline 
$9_6$   &  $ \sqcup_3 (0,0,0),\; 
\sqcup_2   (2 q_1, 0, -2 q_1),\;  \sqcup_2   (0, -2 q_1, 2 q_1),\; 
   \sqcup_2      (-2 q_1, 2 q_1, 0). $ 
\\ \hline 
$9_{10}$   &  $ \sqcup_4 (0,0,0),\; 
 (-q_2 - q_1, 0, q_1 + q_2),\; 
    (q_2, -q_1,  q_1 -q_2),\;  $
\\ \hline 
 & $         (-2 q_1, q_1, q_1),\;  (-q_1 + q_2,  2 q_1 -q_2, -q_1),\; 
    (  q_1-q_2, -2 q_1 + q_2, q_1). $ 
\\ \hline
$9_{11}$   &  $ \sqcup_3 (0,0,0),\;  \sqcup_6 (q_1, 0, -q_1).$
\\ \hline 
$9_{15} $ &  $ \sqcup_3 (0,0,0),\; 
\sqcup_2 (0, -q_1, q_1),\;  \sqcup_2  (-q_1, q_1, 0),\; 
 \sqcup_2  (q_1, 0, -q_1). $  
\\ \hline 
$9_{16} $  &   $ \sqcup_3 (0,0,0),\; 
 \sqcup_2  (-q_1, 0, q_1),\;  \sqcup_2  (q_1, -q_1, 0),\;  
 \sqcup_2  (0, q_1, -q_1).  $ 
\\ \hline 
 $9_{17} $ &  $ \sqcup_3 (0,0,0),\;  \sqcup_6  (-q_1, 0, q_1). $ 
 \\ \hline
$9_{23}$ &  $ \sqcup_3 (0,0,0),\;  
\sqcup_2 (0, -q_1, q_1),\;  \sqcup_2  (-q_1, q_1, 0),\; 
 \sqcup_2  (q_1, 0, -q_1). $
\\ \hline
$9_{24}$ &  $ \sqcup_9 (0,0,0). $ 
\\ \hline
$9_{28}$ &  $ \sqcup_3 (0,0,0),\; 
 \sqcup_3 (-q_1, 0, q_1),\;  \sqcup_3 (0, q_1, -q_1). $ 
\\ \hline
$9_{29}$ &  $ \sqcup_3 (0,0,0),\; 
 (-q_2, -q_1, q_1 + q_2),\;  (q_2, 0, -q_2),\; (0, q_1, -q_1),\;  $
\\ \hline
& $  (-q_1, 0, q_1),\; 
  (q_2,   q_1-q_2, -q_1),\;   (  q_1-q_2,  - q_1+q_2, 0)
. $ 
\\ \hline
 $9_{34}$ &  $ \sqcup_3 (0,0,0),\;  \sqcup_6  (-q_1, 0, q_1). $ 
\\ \hline
$9_{35}$ &  $ \sqcup_3 (0,0,0),\;  \sqcup_8 (-q_1, 0, q_1),\;
 \sqcup_3 (0, q_1, -q_1),\;    (q_1, 0, -q_1),\; (q_1, -q_1, 0),\;  
 (0, -q_2, q_2),\; $
 \\ \hline & $  (0, 2 q_2, -2 q_2),\;
 (-2 q_1, 0, 2 q_1),\; (-q_2,  - 2 q_1-q_2, 2 q_1 + 2 q_2),\; 
  (q_1, q_2, - q_1 -q_2),\; $
 \\ \hline
 & $  (0, q_1 + q_2,  - q_1-q_2),\;
(-q_1, -2 q_2, q_1 + 2 q_2),\;  (-q_1, 2 q_1 - q_2, -q_1 + q_2),\; 
 (q_1, -q_1 + q_2, -q_2),\;  $
 \\ \hline
 & $  (0,   q_1-q_2, -q_1 + q_2),\; 
 (  q_1-q_2, -q_1 + 2 q_2, -q_2). $ 
\\ \hline
$9_{37}$ &  $ \sqcup_9 (0,0,0),\;  \sqcup_5 (0, -q_1, q_1),\;
 \sqcup_5 (q_1, 0, -q_1),\; \sqcup_2 (-q_2, 0, q_2),\;  (-q_2, q_2, 0),$
 \\ \hline
 & $  (-q_2, q_1, -q_1 + q_2),\;  ( q_1+q_2, -q_1, -q_2),\; 
 (q_2, q_1,  - q_1-q_2),\;  (q_2, q_1, - q_1 -q_2), $
 \\ \hline
 & $ (q_1, - 2 q_1 -2 q_2,  q_1+ 2 q_2),\; (q_1, 2 q_2, - q_1 -2 q_2),\;
 (-q_1 + q_2, q_1 - q_2, 0), $
 \\ \hline
 & $ (-q_2, -q_1 + q_2, q_1),\;  (-q_2, -q_1 + q_2, q_1),\; 
  (-q_1 + q_2, q_1 - q_2, 0),\;  (q_2,  - q_1-q_2, q_1). $ 
\\ \hline
 $9_{38}$ &  $ \sqcup_3 (0,0,0),\; 
 (-q_2, -2 q_1, 2 q_1 + q_2),\;   (-q_1 + q_2, q_1, -q_2),\;$
\\ \hline
& $(-q_1, 2 q_1, -q_1),\; (-q_1, -q_1, 2 q_1),\;  (q_2, -q_2, 0),\;
  (-q_2, q_2, 0). $ 
\\ \hline
$9_{40}$ &  $ \sqcup_3 (0,0,0),\; 
\sqcup_3 (-q_1, 0, q_1),\; \sqcup_3(0, q_1, -q_1). $ 
\\ \hline
 $9_{46}$ &  $ \sqcup_{17} (0,0,0),\; 
 \sqcup_3    (q_1, -q_1, 0),\;  \sqcup_3   (0, q_1, -q_1),\;
   (-q_2, q_2, 0),\;    (q_2, -q_2, 0), $
\\ \hline
& $(  q_1+q_2, 0, - q_1 -q_2),\; ( - q_1-q_2, 0,  q_1+ q_2). $ 
\\ \hline
 $9_{47}$ &  $ \sqcup_4 (0,0,0),\; 
\sqcup_{12} (q_1, 0, -q_1),\; \sqcup_3 (0, q_1, -q_1),\;
\sqcup_3  (q_1, -q_1, 0),$
\\ \hline
& $(q_2, -q_2, 0),\;    (-q_2, q_2, 0),\;
 ( - q_1-q_2, 0, q_1 + q_2),\;    (q_1 + q_2, 0, -q_2 - q_1).  $ 
\\ \hline
 $9_{48}$  &  $ \sqcup_3 (0,0,0),\;   \sqcup_4 (0, -q_1, q_1),\;
 \sqcup_3  (q_1, -q_1, 0),\; \sqcup_3  (-q_1, q_1, 0),\; 
\sqcup_2  (0, q_1, -q_1),\; 
\sqcup_2 (q_1, 0, -q_1),\;  
$
\\ \hline
& $  (q_2, -q_2, 0), \;    (-q_2, q_2, 0),
   (-q_1, 2 q_1, -q_1),\;  (-q_1, -q_1, 2 q_1),\;
  (-q_2, -2 q_1, 2 q_1 + q_2),\;    $
\\ \hline
& $  (  q_1+q_2, -q_1, -q_2),\; (-2 q_1 - q_2, q_1,   q_1+q_2),\;  (-q_1 + q_2, q_1, -q_2),\; $
\\ \hline
& $
  (-q_2, -q_1 + q_2, q_1),\; 
  (-q_1 + q_2,  q_1 -q_2, 0). $
\\ \hline
\end{tabular} } \end{center}
\caption{A table of cocycle invariants for twist spun knots, part II}
\label{twistspintable2}
\end{table}


\begin{example} \label{twistspinex} {\rm
We obtained the following calculations by {\it Maple} and {\it Mathematica}. 
Let $X=R_3$,  the coefficient group $\Z^3$ with the wreath
product action (elements of $R_3=\{ 1, 2, 3 \}$ 
acting on $\Z^3$ by transpositions
of factors, $1=(2\ 3)$, $2=(1\ 3)$, and $3=(1\ 2)$). 
 Let $h(i,j, k)={}^T( f_1(i,j,k),  f_2(i,j,k),  f_3(i,j,k)) $
denote a vector valued $3$-cochain, where 
$f_{\ell}(i,j,k)\in \Z$.
Let $q_1, q_2 \in \Z$ be arbitrary elements.
Then the following values, with all the other unspecified ones being zero, 
defined a $3$-cocycle $h$:
\begin{eqnarray*}
f_1(1,2,1) \ =\ f_2(1,2,3) \ =\ f_3(3,1,2) \ =\
- f_1(1,3,1) \ =\ - f_1(1,3,2)   &=& q_1 , \\
f_1(2,1,3)  \ =\ f_1(3,1,2)   \ =\ - f_1(2,3,2)   &=& q_2, \\
f_1(3,2,3)   &=& - q_1 - q_2 
\end{eqnarray*}
With this $3$-cocycle 
the cocycle invariant for the $2$-twist spun knots  for $3$-colorable 
knots in the table are evaluated as listed in
Table~\ref{twistspintable1}  (up to $8$ crossings)
and Table~\ref{twistspintable2} ($9$ crossing knots). 
The notations used in the table for families of vectors 
are similar 
to those  in Example~\ref{r3classicalex}.

 } \end{example} 

\clearpage

\subsection*{Non-invertibility of knotted surfaces} 
The computations of the cocycle invariant imply non-invertibility
of twist-spun knots. Here is a  brief overview on non-invertibility 
of knotted surfaces:
 
\begin{itemize}
\setlength{\itemsep}{-3pt}
\item
\begin{sloppypar}
Fox \cite{Fox61a} presented a  non-invertible knotted sphere
using  knot  modules as follows.
The first homology 
$H_1(\tilde{X})$ of the infinite cyclic cover $\tilde{X}$  
of the complement $X$ of the sphere in $S^4$ 
of Fox's  Example 10 
is   $\Z [ t, t^{-1}]/ (2-t)$   as a $\Lambda=\Z[t,t^{-1}]$-module. 
Fox's Example 11
can be recognized as the same sphere as Example 10 in \cite{Fox61a}
with
its orientation reversed. 
Its Alexander polynomial is $(1-2t)$,
and therefore, not equivalent to Example 10.
\end{sloppypar}

Knot  modules, however,
fail to  detect non-invertibility of
the $2$-twist spun trefoil (whose 
knot module  is 
$\Lambda/(2-t,1-2t)$).

\item
Farber \cite{Far75} 
showed that the $2$-twist spun trefoil  was 
non-invertible  using the Farber-Levine pairing
(see also Hillman \cite{Hill81}).

\item
Ruberman \cite{Rub83}  
used Casson-Gordon invariants  
to prove the same result,
with other new examples of non-invertible knotted spheres. 

\item
\begin{sloppypar}
Neither technique applies
directly to the same knot with 
trivial $1$-handles attached
(in this case the knot is a surface with a higher genus).
Kawauchi~\cite{Kawa86a,Kawa90a}  
has generalized the Farber-Levine pairing to higher genus surfaces,
showing that such a  surface 
is also non-invertible. 
\end{sloppypar}

\item
Gordon~\cite{Gor02*} showed that a large family of knotted spheres
are indeed non-invertible. His extensive lists are:
(1) the $2$-twist spin of a rational knot $K$ is invertible 
if and only if $K$ is amphicheiral;
(2) if $m$, $p$, $q$ are $>1$, then the $m$-twist spin of the $(p,q)$ torus
knot is non-invertible, 
(3) if $m \geq 3$ then the $m$-twist spin of a hyperbolic 
knot $K$ is invertible if and only if $K$ is $(+)$-amphicheiral.
His topological argument uses the fact that the twist spun knots 
are fibered $2$-knots. In particular, the corresponding results are unknown
for surfaces of higher genuses.
\end{itemize}

Then the cocycle invariant provided a diagrammatic method
  of detecting non-invertibility 
of knotted surfaces. In \cite{CJKLS}, it was shown using 
the cocycle invariant that the
$2$-twist spun trefoil is non-invertible.  
Furthermore, as was mentioned in Section~\ref{intro},
all the surfaces that are obtained from  the
$2$-twist spun trefoil
by attaching trivial $1$-handles, called its {\it stabilized} surfaces,
are also non-invertible, since the stabilized surfaces
 have the same cocycle invariant as the original knotted sphere. 
The result was generalized by Satoh~\cite{Sat01b*} to 
an infinite family of twist spins of torus knots. 

The computations given in Example~\ref{twistspinex} 
 can be carried out for the $2$-twist spun knots 
with orientations reversed. 
Specifically, if the orientation of the surface is reversed, then 
the computations change in the following manner. 
First, the face colors are determined by
$$  f(s, j, {\mathcal C}) = ( \cdots ( c(1, k, {\mathcal C}) *
c(s, k, {\mathcal C}) ) * c(s, k-1, {\mathcal C}) )
 *  \cdots )  *
c(s, j, {\mathcal C}) ) \cdots ) .
$$
The target region of $T_1^l(s)$ depicted in Fig.~\ref{triplewt} 
is the top left 
region in the figure, so that 
the sequence of colors that an arc $\gamma$ meets is  
$$
 c(s, n, {\mathcal C})^{-1} c(s, n-1, {\mathcal C})^{-1}  \cdots 
 c(s, j, {\mathcal C})^{-1}, $$
where $\gamma$ does not cross the horizontal sheet in the figure but 
goes through $k$-th through $j$-th vertical sheets from left. 
Thus we obtain
$$
  B(  {\mathcal C}, T_1^l(s))
= \left\{ 
\begin{array}{ll}
 c(s, n, {\mathcal C})^{-1}  \cdots  c(s, j, {\mathcal C})^{-1} 
\kappa_{f(s, j+2, {\mathcal C}),\
c(s+1, j+1,  {\mathcal C}),\  c(s, j+1,  {\mathcal C})} & \quad
\mbox{\rm if}\quad  \epsilon (T_1^l(s)) > 0 
\\
-( 
 c(s, n, {\mathcal C} )^{-1}  \cdots  c(s, j, {\mathcal C})^{-1}  )
\kappa_{f(s, j+2, {\mathcal C}),\ c(s, j+1, {\mathcal C}),\ 
c(s, j, {\mathcal C})
} 
 & \quad
\mbox{\rm if}\quad \epsilon (T_1^l(s)) < 0 ,
\end{array}
\right.
$$ 
$$ 
  B(  {\mathcal C}, T_1^r(s))
= \left\{ 
\begin{array}{ll}
-(
 c(s, n, {\mathcal C})^{-1}  \cdots  c(s, j, {\mathcal C})^{-1} )
\kappa_{c(s+1, j+1, {\mathcal C})
\ c(s, j+1, {\mathcal C}) 
\ c(1, k, {\mathcal C})}  \quad &
\mbox{\rm if} \quad \epsilon (T_1^r(s)) < 0 
\\
 c(s, n, {\mathcal C})^{-1}   \cdots  c(s, j, {\mathcal C})^{-1}  
\kappa_{c(s, j+1, {\mathcal C})
\ c(s, j, {\mathcal C}) 
\
c(1, k, {\mathcal C})} 
  \quad &
\mbox{\rm if} \quad \epsilon (T_1^r(s)) > 0 . 
\end{array}
\right.
$$

The computational results are presented in Tables~\ref{twistspintable3} and 
~\ref{twistspintable4}.
By comparing the computational results, we conclude that such
$2$-twist spun knots $\tau^2(K)$ are non-invertible, for those knots that
give rise to distinct values for the cocycle invariants, 
as well as  
all of their stabilized surfaces. 
This list (with stabilized surfaces) of non-invertible 
surfaces has 
not been obtained  by 
any other method.

Furthermore, it is easily seen that if  the contribution 
to the invariant for a coloring ${\cal C}$ is a vector $(a,b,c) \in \Z^3$ for
the $2$-twist spun  $\tau^2(K)$ of a classical knot $K$,
then the contribution for the corresponding coloring  ${\cal C}$
is the vector  $k (a,b,c)$ for the $2k$-twist spun   $\tau^{2k}(K)$
of $K$. Hence non-invertibility determined by this invariant for 
$2$-twist spuns can be applied to $2k$-twist spuns for all
positive integer $k$ as well. 

We summarize the result in the following theorem.

\begin{theorem}
For any positive integer $k$, 
the $2k$-twist spun  of all the 
$3$-colorable  knots in the table up to $9$ crossings  excluding $8_{20}$, 
as well as their stabilized surfaces of any genus, 
are non-invertible.
\end{theorem}

The cocycle invariant in Example~\ref{twistspinex} 
fails to detect non-invertibility of  $8_{20}$ and 
we obtain no conclusion.

The $2$-twist spins of $24$ of these $3$-colorable knots 
can be distinguished from their inverses  
by the invariant of \cite{CJKLS}. Their $6k$-twist spins,
however, as well as the $2k$-twist spins of the following list
and their stabilizations,
can be distinguished from their inverses only by the current 
invariant: 
$6_1$, $8_{10}$,  $8_{11}$, $8_{18}$, 
$9_1$, $9_6$, $9_{23}$, $9_{24}$,  $9_{37}$, $9_{46}$.

Furthermore, the original invariant for $R_3$ with trivial action has 
non-trivial cohomology only with $\Z_3$, so that the  invariant 
is trivial for $6k$-twist spuns, and in particular, 
unable to detect non-invertibility 
in this case. The current invariant with non-trivial action,
having a free abelian coefficient group $\Z^3$, is non-trivial for 
any $2k$-twist spuns, if the invariant 
is non-trivial for the $2$-twist spun.

\subsection*{Acknowledgement} We thank 
 Pat Gilmer, Seiichi Kamada, Dan Silver, and Susan Williams
for valuable comments and discussions.

\begin{table}[h] 
\begin{center}
{\begin{tabular}{||l||l||}\hline \hline
Knot $K$ & $\Phi_{\kappa}({\rm Tw}^2 (K) )$
\\ \hline \hline 
$3_1$ & $\sqcup_3 (0,0,0),\; 
\sqcup_2 (q_1, 0, -q_1),\; \sqcup_2 (0, -q_1, q_1),\; 
\sqcup_2  (-q_1, q_1, 0). $ 
\\ \hline 
$6_1$ & $\sqcup_3 (0,0,0),\; 
( 0,q_1, -q_1), \; ( q_1, -q_1, 0), \; ( 0,-q_2, q_2),\;  $
\\ \hline 
 & $(q_1, q_2, -q_1-q_2),\; (q_1,  -q_1+q_2, -q_2),\; (0, q_1-q_2,-q_1+ q_2).$  
\\ \hline 
$7_4$  & $ \sqcup_3 (0,0,0),\; 
  (0,   q_1+q_2, - q_1 -q_2),\;    (q_1, -q_2, -q_1 + q_2),  $
\\ \hline 
 & $ (2 q_1, -q_1, -q_1),\;  (q_1, -2 q_1 + q_2,   q_1-q_2),\; 
       (-q_1,  2 q_1-q_2 , -q_1 + q_2). $ 
\\ \hline 
$7_7$    &  $ \sqcup_3 (0,0,0),\;
 (q_1,   q_1+q_2, - 2 q_1 -q_2),\; 
        (0,  - 2 q_1-q_2, 2 q_1 + q_2),\;  (0, -q_1, q_1),
$ \\ \hline 
 & $   (q_1, 0, -q_1),\;  (q_1, q_2, - q_1 -q_2),\; 
       (0, - q_1 -q_2,   q_1+q_2).  $ 
\\ \hline 
$8_5$  &  $ \sqcup_3 (0,0,0),\;  
\sqcup_3 (2 q_1, -2 q_1, 0),\; \sqcup_3  (0, 2 q_1, -2 q_1).$ 
\\ \hline 
$8_{10}$ &  $ \sqcup_3 (0,0,0),\;
\sqcup_2 (-q_1, 0, q_1),\; (0, q_1, -q_1),\; (q_1, -q_1, 0),\;
$ \\ \hline 
 & $  (-2 q_1, 2 q_1, 0),\;       (0, -2 q_1, 2 q_1). $ 
\\ \hline
$8_{11}$   &  $ \sqcup_3 (0,0,0),\; \sqcup_2  (0, -q_2, q_2),\;
 (q_1, -2 q_1, q_1),\; (q_1, q_1, -2 q_1),
$ \\ \hline 
 & $(-q_1,   q_1-q_2, q_2),\; 
  (2 q_1, q_2, - 2 q_1 -q_2).  $
\\ \hline 
$8_{15}$   &  $ \sqcup_3 (0,0,0),\;  
 \sqcup_2  (q_1, 0, -q_1), \;   \sqcup_2 (-q_1, q_1, 0),\;  
 \sqcup_2  (0, -q_1, q_1). $  
\\ \hline 
$8_{18}$    &  $ \sqcup_9 (0,0,0),\;
 \sqcup_3 (0,- q_1 -q_2 ,  q_1 + q_2 ),\;
 \sqcup_2  (q_1, 0, -q_1), \;
\sqcup_2    (0, -q_1, q_1),  \;  (-q_1, 2 q_1, -q_1), \;   
$
\\ \hline 
    &  $  \sqcup_2    (q_1,  q_1+ q_2, - 2 q_1 -q_2 ), \;
  \sqcup_2   (0, - 2 q_1 -q_2, 2 q_1 + q_2), \;
 (q_1, -q_1 + q_2, -q_2), \; $
\\ \hline 
    &  $ (q_1, q_2, -q_2 - q_1), \;
      (2 q_1, -q_1 + q_2, - q_1 -q_2), \;  
       (q_1, q_1 - q_2, -2 q_1 + q_2), \;
     (q_1, q_2,  - q_1-q_2).$
\\ \hline
$8_{19}$    &   $ \sqcup_5 (0,0,0),\;
 (-q_1, 2 q_1, -q_1),\; (2 q_1, -q_1, -q_1),\;
        (q_1, -2 q_1, q_1),\; (q_1, q_1, -2 q_1).
 $ 
\\ \hline 
$8_{20} $   &  $ \sqcup_9 (0,0,0).  $  
\\ \hline 
$8_{21}$    &   $ \sqcup_5 (0,0,0)\;  (q_1, q_1, -2 q_1), \;     
 (2 q_1, -q_1, -q_1),  \; (-q_1, 2 q_1, -q_1),\; (q_1, -2 q_1, q_1).
 $  
\\ \hline 
\end{tabular} } \end{center}
\caption{A table of cocycle invariants
with orientations reversed, part I}
\label{twistspintable3}
\end{table}

\begin{table}[h] 
\begin{center}
{\begin{tabular}{||l||l||}\hline \hline
Knot $K$ & $\Phi_{\kappa}({\rm Tw}^2 (K) )$
\\ \hline \hline 
$9_1$   &    $ \sqcup_3 (0,0,0),\; 
\sqcup_2 (3 q_1, 0, -3 q_1),\; \sqcup_2  (-3 q_1, 3 q_1, 0),\;
\sqcup_2  (0, -3 q_1, 3 q_1).
 $ 
\\ \hline
$9_{2}$    &   $ \sqcup_3 (0,0,0),\; 
 (-q_2, 0, q_2),\;    (q_1, 0, -q_1),\;  (2 q_1, -2 q_1, 0),\; 
    (q_2,  - q_1 -q_2, q_1),\; $
\\ \hline  &
$ ( - q_1-q_2,  q_1+q_2, 0),\; (q_1+q_2, 2 q_1, -3 q_1 - q_2). 
 $  
\\ \hline
$9_4$     &   $ \sqcup_3 (0,0,0),\;  (0, -2 q_1, 2 q_1),\;   (0, q_1, -q_1), $
\\ \hline  &
$ (-2 q_1,  q_1 -q_2,  q_1+q_2 ),\;
      (q_1, q_2,  - q_1-q_2),\;
   (-q_1, -q_2,  q_1+q_2 ),\;
      (-q_1, q_2,  q_1-q_2 ).
 
 $  
\\ \hline 
$9_6$   &  $ \sqcup_3 (0,0,0),\;  (0, -3 q_1, 3 q_1),\; (3 q_1, 0, -3 q_1),\;
 (q_1, 2 q_1 + q_2, - 3 q_1 -q_2),$
\\ \hline & $ (q_1,  - 3 q_1-q_2, 2 q_1 + q_2),\;
     (2 q_1, -q_1 + q_2,  - q_1-q_2),\;
     (-q_1,  - q_1 -q_2, 2 q_1 + q_2). 
 $ 
\\ \hline 
$9_{10}$   &  $ \sqcup_3 (0,0,0),\; 
 (2 q_1, q_2, -q_2 - 2 q_1),\; 
      (-q_1, -q_2 + q_1, q_2),\; $
\\ \hline &
$ (q_1, q_1, -2 q_1),
        (q_1, -2 q_1, q_1), (0, q_2, -q_2), (0, -q_2, q_2).
 $
\\ \hline
$9_{11}$   &  $ \sqcup_3 (0,0,0),\; 
 \sqcup_3  (q_1, -q_1, 0),\;  \sqcup_3  (0, q_1, -q_1).
$
\\ \hline 
$9_{15} $ &  $ \sqcup_3 (0,0,0),\; 
 \sqcup_2 (0, q_1, -q_1),\;  \sqcup_2 (2 q_1, 0, -2 q_1),\; 
 \sqcup_2        (q_1, -q_1, 0). 
       $  
\\ \hline 
$9_{16} $  &   $ \sqcup_3 (0,0,0),\; 
 \sqcup_2  (q_1, 0, -q_1),\;  \sqcup_2   (-q_1, q_1, 0), \; 
  \sqcup_2    (0, -q_1, q_1).
  $ 
\\ \hline 
 $9_{17} $ &  $ \sqcup_9 (0,0,0).  
 $ 
 \\ \hline
$9_{23}$ &  $ \sqcup_5 (0,0,0),\;  
 (q_1, q_1, -2 q_1),   \;    (2 q_1, -q_1, -q_1),\;  (q_1, -2 q_1, q_1),\; 
   (-q_1, 2 q_1, -q_1). 
 $
\\ \hline
$9_{24}$ &  $ \sqcup_3 (0,0,0),\;     \sqcup_2    (-q_1, 0, q_1),\;
 (q_1, -q_1, 0),\;    (0, q_1, -q_1),\; (0, -2 q_1, 2 q_1),\;    (-2 q_1, 2 q_1, 0). $ 
\\ \hline
$9_{28}$ &  $ \sqcup_6 (0,0,0),\; 
 (-q_1, 2 q_1, -q_1),\; 
       (2 q_1, -q_1, -q_1), \; 
    (-q_1, -q_1, 2 q_1).
 $ 
\\ \hline
$9_{29}$ &  $ \sqcup_3 (0,0,0),\;   (0, q_1, -q_1),\;
 (0, -q_2, q_2),\;     (q_1, -q_1, 0),$
\\ \hline & 
$  (q_1, q_2, -q_2 - q_1), \;    (q_1, q_2 - q_1, -q_2), \;  
(0, -q_2 + q_1, q_2 - q_1).
 $
\\ \hline
 $9_{34}$ &  $ \sqcup_9 (0,0,0).  
 $ 
\\ \hline
$9_{35}$ &  $ \sqcup_4 (0,0,0),\; \sqcup_3  (0, -q_1, q_1),\;
\sqcup_3  (-q_1, q_1, 0),\;
\sqcup_3 (0, 2 q_1, -2 q_1),\; 
 \sqcup_3 (2 q_1, -2 q_1, 0),\;$
 \\ \hline & 
$  (q_1, q_1, -2 q_1),\;    (q_1, -2 q_1, q_1), \; (2 q_1, -q_1, -q_1),\;
 (-q_2, -q_1,   q_1+q_2),\;$
 \\ \hline & 
$  (q_1, 2 q_2,  - q_1-2 q_2),\;  (q_1 - q_2,  - q_1+ q_2, 0),\;
  (  q_1+q_2,  - q_1-q_2, 0),\; 
        (  q_1+q_2,   q_1+q_2, -2 q_1 - 2 q_2), \;$
 \\ \hline & 
$    (q_1, -2 q_1 - 2 q_2,   q_1+2 q_2),\;
     (  2 q_1+q_2, -q_1, - q_1 -q_2), \;
      ( - q_1+q_2,   2 q_1-2 q_2,  - q_1+q_2).
 $
 \\ \hline
$9_{37}$ &  $ \sqcup_3 (0,0,0),\;
\sqcup_2 (0, q_1, -q_1),\; \sqcup_2   (q_1, -q_1, 0),\;
\sqcup_2  (0, -q_2, q_2),\; \sqcup_3 (0, -q_2 + q_1, -q_1 + q_2), $
 \\ \hline &
$ (0, -q_1, q_1), \; (-q_1, q_1, 0), \; (0, q_2, -q_2),
 \;(0, 2 q_2, -2 q_2),\; 
(0, -q_1 + q_2, -q_2 + q_1),\;$
 \\ \hline &
$    (q_1, -q_2, -q_1 + q_2),\;
   (q_1, q_2, - q_1 -q_2),\;
     (2 q_1 + q_2, 0, -2 q_1 - q_2),\;
    (q_1, q_2,  - q_1-q_2),\;$
 \\ \hline &
$ 
      ( - q_1-q_2, -q_1, 2 q_1 + q_2),\;  (-q_1, -2 q_2,  q_1+ 2 q_2),\;
    (q_1, -q_1 + q_2, -q_2),\;
    ( - q_1-q_2, q_1 + q_2, 0), \;$
 \\ \hline &
$ 
    (q_1 + q_2, -q_2, -q_1),\;  (q_1, -q_1 + q_2, -q_2).
 $
 \\ \hline
 $9_{38}$ &  $ \sqcup_4 (0,0,0),\; 
 (0,  q_1+q_2, - q_1 -q_2),\; 
   (q_1, -q_2, -q_1 + q_2),\; 
    (2 q_1, -q_1, -q_1), \;$
 \\ \hline &
$  (q_1, -2 q_1 + q_2, q_1 - q_2),\; 
     (-q_1, 2 q_1 - q_2, -q_1 + q_2). 
 $
\\ \hline
$9_{40}$ &  $ \sqcup_6 (0,0,0),\; 
 (-q_1, 2 q_1, -q_1),\; 
    (2 q_1, -q_1, -q_1), \; 
      (-q_1, -q_1, 2 q_1).
 $ 
\\ \hline
 $9_{46}$ &  $ \sqcup_{15} (0,0,0),\; 
\sqcup_2 (q_1, -q_1, 0),\; \sqcup_2 (0, -q_1, q_1),\;
\sqcup_2  (-q_1, 0, q_1),\;  (-q_1, q_1, 0), \;  (0, q_1, -q_1),\;  (0, q_2, -q_2),\;$
\\ \hline &  
$  (-q_1, -q_2,   q_1+ q_2), \; 
  (-q_1,  q_1 -q_2, q_2),  \; 
      (0, - q_1+ q_2,  q_1 -q_2). 
  $
\\ \hline
 $9_{47}$ &  $ \sqcup_{11} (0,0,0),\;
\sqcup_3 (0, -q_2, q_2),\; \sqcup_2  (0, q_1, -q_1),\;
\sqcup_2  (q_1, -q_1, 0), \; 
 \sqcup_2   (0,   q_1-q_2,  - q_1+ q_2),\;  (0, q_2, -q_2),\; $
\\ \hline &
$\sqcup_2  (q_1, q_2, - q_1 -q_2 ),\;
 \sqcup_2    (q_1,  - q_1+q_2, -q_2),\;
  (0,  q_1+ q_2,  - q_1-q_2),\;
     (0,  - q_1-q_2,   q_1+q_2).
$
\\ \hline
 $9_{48}$  &  $ \sqcup_3 (0,0,0),\;   
  \sqcup_6 (0, -q_1, q_1),\;   
  \sqcup_4 (q_1, 0, -q_1),\;
  \sqcup_4 (-q_1, q_1, 0),$
\\ \hline &
$ (0, q_1, -q_1), \;   
  (q_1, -q_1, 0),\;   
  (0, q_2, -q_2),\;   
   (0, -q_2, q_2),$
\\ \hline &
$
 (q_1, q_2, - q_1 -q_2),\;   
 (q_1, - q_1 -q_2, q_2),\;   
  (q_1,  - q_1+ q_2, -q_2),$
\\ \hline &
$
 (0,   2 q_1-q_2,  - 2 q_1+ q_2),\;   
  (2 q_1,  - 2 q_1+ q_2, -q_2),\;   
   (0,  q_1 -q_2, - q_1+ q_2).  
$
\\ \hline
\end{tabular} } \end{center}
\caption{A table of cocycle invariants with orientation reversed, part II}
\label{twistspintable4}
\end{table}

\end{document}